\documentclass[reqno, 11pt]{amsart}  
\usepackage[left=1.15in,right=1.15in,top=1in,bottom=1in]{geometry}

\usepackage{amsmath}
\usepackage{amssymb}
\usepackage{hyperref}
\hypersetup{colorlinks=true}
\usepackage{graphicx}
\graphicspath{~/Desktop/first-order hyperbolic/}
\numberwithin{equation}{section}

\newtheorem{prop}{Proposition}[section]
\newtheorem{lemma}[prop]{Lemma}
\newtheorem{thm}[prop]{Theorem}
\newtheorem{cor}[prop]{Corollary}

\newtheorem*{helton_vinnikov_thm}{Helton-Vinnikov Theorem}

\newcommand{\shrink}[1]{ {\scriptstyle {\textstyle {#1} } } }
\newcommand{\smfrac}[2]{ \shrink{ \frac{#1}{#2} } }

\newcommand{\lin}{\langle}
\newcommand{\rin}{\rangle}
\newcommand{\Sym}{\mathbb{S}^{n}}
\newcommand{\Symp}{\Sym_{\plus}} 
 
\newcommand{\realsnp}{\mathbb{R}^n_{\plus}}

\newcommand{\plus}{{\scriptscriptstyle +}}

\renewcommand{\int}{\mathrm{int}}

\newcommand{\cp}{\mathrm{CP}} 

\newcommand{\lambdamin}{\lambda_{\mathrm{min}}} 
\newcommand{\lambdamax}{\lambda_{\mathrm{max}}} 
\newcommand{\grad}{\nabla} 

\newcommand{\diamval}{\mathrm{diam}_{z}  }

\newcommand{\hatfmu}{\hat{f}_{\mu}} 

\begin{document}

\newpage 
$ \textrm{~} $ \quad \vspace{-3mm}

\title[Accelerated Methods 
  for Hyperbolic Programming]{Accelerated First-Order Methods \\
  for Hyperbolic Programming}

\begin{abstract}
We develop a framework for applying accelerated methods to general hyperbolic programming, including linear, second-order cone, and semidefinite programming as special cases. The approach replaces a hyperbolic program with a convex optimization problem whose smooth objective function is explicit, and for which the only constraints are linear equations (one more linear equation than for the original problem).  Virtually any first-order method can be applied. An iteration bound for a representative accelerated method is derived.
\end{abstract} \vspace{-3mm}

\author[J. Renegar]{James Renegar  }
\address{School of Operations Research and Information Engineering,
 Cornell University, Ithaca, NY, U.S.}
 
 \thanks{Research supported in part by NSF CCF Award \#1552518.}

\maketitle

\vspace{-7mm}

\section{{\bf  Introduction}}  \label{sect.a}
In \cite{renegar2016efficient}, we showed that given a general  conic optimization problem for which a strictly feasible point is known, there is a  simple transformation to an equivalent convex optimization problem which has the same number of variables, has only linear equations as constraints (one more linear equation than the original problem), and has Lipschitz-continuous objective function defined on the whole space. Virtually any subgradient method can be applied to solve the equivalent problem, the cost per iteration dominated by computation of a subgradient and its orthogonal projection onto a subspace (the same subspace at every iteration, a situation for which preprocessing is effective). Iteration bounds for a representative subgradient method were established. \vspace{1mm}

The iteration bounds, however, are of order $ 1/\epsilon^2 $, as is to be expected for subgradient methods when the goal is to compute a point whose objective value is within $ \epsilon $ of optimal (\cite{nemirovsky1983problem}).  By contrast, a decade ago, Nesterov \cite{nesterov2005smooth}  obtained $ O(1/ \epsilon) $ iteration bounds for an accelerated first-order method applied to various minimization problems in which the objective is the sum of a smooth convex function and a nonsmooth simple convex function, where by ``simple'' is meant that the function is well-approximated by an implicitly-defined smooth function (a Moreau proximal smoothing) for which the gradient at a point can be computed by solving a simple optimization problem  -- an optimization problem so simple that the gradient can be expressed in closed form, that is, made explicit.     \vspace{1mm}

This work of Nesterov motivated many researchers, leading to a great number of related publications. An elegant formalization and insightful overview was provided by Beck and Teboulle \cite{beck2012smoothing}, indicating the generality of the approach, at least in theory.  The limitation of the approach is that for the implicitly-defined smooth function, the subproblem for computing the gradient at a point indeed has to be simple, else the cost of (approximately) computing the gradient can greatly exceed all other costs in an iteration.  This requirement of simplicity significantly restricts the scope of optimization problems to which the approach can be applied in practice, although the approach has been very consequential for some particularly important classes of high-profile optimization problems, such as in compressive sensing (c.f., \cite{becker2011nesta, becker2011templates, foucart2013mathematical}). \vspace{1mm}

Here we extend our investigation in \cite{renegar2016efficient}, now applying accelerated methods rather than subgradient methods.  However, whereas \cite{renegar2016efficient} considered  {\em all} convex conic optimization problems, in this paper attention is restricted to hyperbolic programming, still an extremely-general class of conic optimization problems, albeit a class for which the algebraic structure of the underlying cones naturally provides an explicit ``smoothing'' of the objective function for the equivalent optimization problem.  Once the objective function is smoothed, as will become clear, virtually any first-order method can be applied.  
\vspace{1mm}

Semidefinite programs are hyperbolic programs, and hence so are linear programs and second-order cone programs.  It remains unknown, however, whether hyperbolic programming is more general than semidefinite programming, that is, it remains unknown whether every hyperbolic program might be recast as a semidefinite program, possibly involving more variables.  A primary reason we focus on hyperbolic programming is to make our results be flexible for practice.  Our general recipe for computing gradients of the smoothed objective function, for example, immediately specializes to a wide variety of settings, without the need, say, to first convert problems to equivalent semidefinite programs, and without the need of introducing additional variables.  \vspace{1mm}

In  \S\ref{sect.ba}, we provide a brief overview of the framework and ``basic theory'' introduced in \cite[\S 2]{renegar2016efficient}, all of which is elementary.   In \S\ref{sect.bb}, we give a concise introduction to hyperbolic programming, a task easily accomplished by relying on the reader's grasp of the basics for semidefinite programming.  \vspace{1mm}

Results regarding ``smoothing'' are presented in \S \ref{sect.c}  (lengthy proofs are deferred to appendices).  Besides presenting desirable theoretical properties of the explicit smoothing, it is here that we provide computational expressions for gradients, making evident the practicality of the smoothing especially when a hyperbolic program has a ``hyperbolic polynomial'' which factors as a product of polynomials of low degree.  The results are demonstrated by showing for a general second-order cone program, computing gradients for the smoothed function \vspace{1mm} is easy.

In \S\ref{sect.d}  we consider a representative accelerated method, due to Nesterov \cite{nesterov2013gradient}, establishing an $ O(1/ \epsilon) $-iteration bound when the method is applied to the smoothing of a general hyperbolic program. Additionally, we give considerable attention to the constant hidden in the big-$ O $, showing that if the accelerated method is instead applied to a carefully-chosen {\em  sequence} of  smoothed versions of the hyperbolic program, significantly greater efficiency can be obtained, that is, far fewer gradient computations can suffice for solving the hyperbolic program to within desired accuracy.
\vspace{1mm}

  These results significantly extend and strengthen ones announced in \cite{renegar2014efficient}  for semidefinite programming.  The present paper and \cite{renegar2016efficient}  should together be considered as the final, refined version of \cite{renegar2014efficient}.\vspace{1mm}
  
An especially-notable earlier work with a similar $ O(1/ \epsilon) $-iteration bound for the special case of semidefinite programming is due to Lu, Nemirovski and Monteiro \cite{lu2007large}  (rather than a framework in which various accelerated methods can be applied, their focus is on a particular accelerated method aimed at solving large-scale problems). Also deserving of special mention is work of Lan, Lu and Monteiro \cite{lan2011primal}, who obtain $ O(1/\epsilon)$-iteration bounds for primal-dual convex conic optimization problems (although for the cost of  individual iterations to be low, projection onto the underlying cones has to be inexpensive). \vspace{1mm}

Closely related to the present paper and \cite{renegar2016efficient}  is recent significant work of Freund and Lu \cite{freund2015new}, who develop first-order methods for convex optimization problems in which the objective function has a known lower bound. They consider general non-smooth, smooth and smoothed settings (although for the latter, they assume the smoothing is given, whereas the relevance of the present paper rests largely in contributing a particular smoothing). The perspective they develop provides an important juxtaposition to our conic-oriented presentation.

\section{{\bf  Preliminaries}}  \label{sect.b}

\subsection{{\bf Basic Theory}} \label{sect.ba}

Here we provide a brief overview of the general framework introduced in \cite[\S 2]{renegar2016efficient} and record relevant results proven there. All of this is elementary. \vspace{1mm}

Let $ {\mathcal E} $ be a finite-dimensional real Euclidean space.  Let $ {\mathcal K}  \subset {\mathcal E} $ be a proper, closed, convex cone with non-empty interior. \vspace{1mm}

Fix a vector $ e \in \int(  {\mathcal K} ) $ (interior).  We refer to $ e $ as the ``distinguished direction.''   For each $ x \in {\mathcal E} $, let 
\[   \lambda_{\min}(x) :=  \inf \{ \lambda: x - \lambda \, e \notin {\mathcal K}  \} \; , \]
that is, the scalar $ \lambda $ for which $ x - \lambda e $ lies in the boundary of $ {\mathcal K}  $.  (Existence and uniqueness of $ \lambda_{\min}(x) $ follows from $ e \in \int({\mathcal K} ) \neq {\mathcal E} $ and convexity of $ {\mathcal K}  $.)  \vspace{1mm}

If, for example, $ {\mathcal E} = \Sym $ ($ n \times n $ symmetric matrices), $ {\mathcal K}  = \Symp $ (cone of positive semidefinite matrices), and $ e = I $ (the identity), then $ \lambda_{\min}(X) $ is the minimum eigenvalue of $ X $. \vspace{1mm}

On the other hand, if $ {\mathcal K}  = \realsnp $ (non-negative orthant) and  $ e $ is a vector with all positive coordinates, then $ \lambda_{\min}(x) = \min_j x_j/e_j $ for $ x \in \mathbb{R}^n $. Clearly, the value of $ \lambda_{\min}(x) $ depends on the distinguished direction $ e $ (a fact the reader should keep in mind since the notation does not reflect the dependence).  \vspace{1mm}

Obviously, $ {\mathcal K}  = \{ x: \lambda_{\min}(x) \geq 0 \} $ and $ \int({\mathcal K} ) = \{ x: \lambda_{\min}(x) > 0 \} $. Also,
\begin{equation}  \label{eqn.ba}
  \lambda_{\min}(sx + te) = s \, \lambda_{\min}(x) + t \quad \textrm{for all $ x \in {\mathcal E}$  and scalars $ s \geq 0 $,  $ t $} \; .  
  \end{equation}

Let 
\[ \bar{{\mathcal B}} := \{ v \in {\mathcal E}: e + v, e - v \in {\mathcal K}  \} \; , \]
a closed, centrally-symmetric, convex set with nonempty interior. 
Define a seminorm\footnote{Recall that a seminorm $ \| \, \, \| $ satisfies $ \| tv \| = |t| \, \| v \| $ and $ \| u + v \| \leq \| u \| + \| v \| $, but unlike a norm, is allowed to satisfy $ \| v \| = 0 $ for $ v \neq 0 $.}  on $ {\mathcal E} $ according to 
\begin{equation}  \label{eqn.bb} 
  \| u \|_{\infty} := \min \{ t \geq 0 : u = tv \textrm{ for some $ v \in \bar{{\mathcal B}} $} \} \; . 
  \end{equation} 
Let $ \bar{B}_{\infty}(x,r) $ denote the closed ball centered at $ x $ and of radius $ r $. Clearly, $ \bar{B}_{\infty}(0,1) = \bar{{\mathcal B}} $, and $ \bar{B}_{\infty}(e,1) $ is the largest subset of $ {\mathcal K}  $ that has symmetry point $ e $, i.e., for each $ v $, either both points $ e + v $ and $ e - v$ are in the set, or neither point is in the set.      
\vspace{1mm}

It is straightforward to show $ \| \, \, \|_{\infty} $ is a norm if and only if $ {\mathcal K}  $ is pointed (i.e., contains no subspace other than $ \{ \vec{0}\} $). \vspace{1mm}

\begin{prop} (\cite[Prop 2.1]{renegar2016efficient})  \label{prop.ba}
  The function $ x \mapsto \lambda_{\min}(x) $ is concave and Lipschitz continuous with Lipschitz constant $ 1 $: 
\[    | \lambda_{\min}(x) - \lambda_{\min}(y) | \leq \| x - y \|_{\infty} \quad \textrm{for all $ x,y \in {\mathcal E} $} \; .   \]  
\end{prop} 
 \vspace{3mm}

Assume the Euclidean space $ {\mathcal E} $ is endowed with an inner product written $ u \cdot v $.  Let  $ \mathrm{Affine} \subseteq {\mathcal E} $ be an affine space, i.e., the translate of a subspace. For fixed $ c \in {\mathcal E} $, consider the conic program
\[  
  \left. \begin{array}{rl}
\inf & c \cdot x  \\
\textrm{s.t.} & x \in \mathrm{Affine}  \\
  & x \in {\mathcal K}   \end{array} \right\} \cp  \]  
Let $ z^* $ denote the optimal value. \vspace{1mm}

Assume $ c $ is not orthogonal to the subspace of which $ \mathrm{Affine} $ is a translate, since otherwise all feasible points are optimal.  This assumption implies that all optimal solutions for CP lie in the boundary of $ {\mathcal K}  $. \vspace{1mm}

Assume $ \mathrm{Affine} \cap \int({\mathcal K} ) $ -- the set of strictly feasible points -- is nonempty.  Fix a strictly feasible point, $ e $.  The point $ e $ serves as the distinguished direction. \vspace{1mm}

For scalars $ z \in \mathbb{R}  $, we introduce the affine space
\[ \mathrm{Affine}_{z} := \{ x \in \mathrm{Affine}:   c \cdot x = z \} \; . \]
Presently we show that for any choice of $ z $   satisfying $ z < c \cdot e \; ,  $ CP can be easily transformed into an equivalent optimization problem in which the only constraint is $ x \in \mathrm{Affine}_{z} \; . $ We need a simple assumption, namely that CP has bounded optimal value, which leads to the following simple observation.   \vspace{1mm}

\begin{lemma}  \label{lem.bb} 
Assume $ \mathrm{CP} $ has bounded optimal value.  

$ \textrm{~} $ \qquad \qquad  \qquad  If $ x \in \mathrm{Affine} $ satisfies $ c \cdot x < c \cdot e $, then $ \lambdamin(x) < 1 \; . $ 
\end{lemma}
\noindent {\bf Proof:} If $ \lambdamin(x) \geq 1 $, then by (\ref{eqn.ba}), $ \lambdamin(e + t(x-e)) \geq 1 $ for all $ t \geq 0 $, and hence $ x + t(x-e) $ is feasible for all $ t \geq 0 $. Since $ \lim_{t \rightarrow \infty} c \cdot (e + t(x-e)) = -\infty $, the lemma is immediate.   \hfill $ \Box $
 \vspace{3mm}
 
For  $ x \in {\mathcal E}  $ satisfying $ \lambda_{\min}(x) < 1 $, let $ \pi(x) $ denote the point where the half-line beginning at $ e $ in direction $ x - e $ intersects the boundary of $ {\mathcal K}  $:
\[  \pi(x) :=  e + \smfrac{1}{1 - \lambda_{\min}(x)} (x - e)    \]
(to verify correctness of the expression, observe (\ref{eqn.ba}) implies $ \lambda_{\min}(\pi(x)) = 0 $). 
We refer to $ \pi(x) $ as the ``radial projection'' (from $ e $) of $ x $ to the boundary of the \vspace{1mm}  feasible region. 

The following theorem has a simple proof, and yet the theorem is the foundation of both \cite{renegar2016efficient} and the present paper. \vspace{1mm}

\begin{thm} (\cite[Thm 2.3]{renegar2016efficient}) \label{thm.bc}
Let $ z $ be any value satisfying \, $   z < c \cdot e \; . $  If $ x^* $ solves
\begin{equation}  \label{eqn.bc}
  \begin{array}{rl}
   \sup & \lambdamin(x) \\
 \mathrm{s.t.} &  x \in \mathrm{Affine}_{z} \; ,  \end{array} 
    \end{equation} 
then $ \pi( x^* ) $ is optimal for $ \mathrm{CP} $. Conversely, if $ \pi^*  $ is optimal for $ \mathrm{CP} $, then $ x^* :=   e + \frac{c \cdot e - z }{c \cdot e - z^* } ( \pi^* - e) $ is optimal for (\ref{eqn.bc}), and $ \pi^*  = \pi( x^* ) $.  
\end{thm} \vspace{1mm}

Under the assumption that a strictly feasible point $ e $ is known, CP has been transformed into an equivalent linearly-constrained maximization problem with concave, Lipschitz-continuous objective function. Virtually any subgradient method --  rather, {\em  supgradient} method -- can be applied to this problem, the main cost per iteration being in computing a supgradient and projecting it onto the subspace $ {\mathcal L}  $ of which the affine space $ \mathrm{Affine}_{z} $ is a translate.  This approach was investigated in \cite{renegar2016efficient}.  In the present paper, where the focus is on accelerated gradient methods, the objective function of the equivalent problem (\ref{eqn.bc})  is smoothed.  Before discussing smoothing, however, we complete preliminaries.  \vspace{1mm}

{\em  Assume, henceforth, that CP has at least one optimal solution, and that $ z $ is a fixed scalar satisfying \, $ z < c \cdot e $.}   Then the equivalent problem (\ref{eqn.bc}) has at least one optimal solution.     Let $ \lambda_{\min,z}^*  $ denote the optimal value for the equivalent problem, and recall
$ z^* $ denotes the optimal value of CP. A useful characterization of $ \lambda_{\min,z}^* $ is easily provided. \vspace{1mm}

\begin{lemma} (\cite[Lemma 2.4]{renegar2016efficient}) \label{lem.bd}
 \[   \lambda_{\min,z}^* = \frac{z - z^* }{c \cdot e - 
 z^* } \; . \] 
  \end{lemma}
\vspace{2mm}

We focus on the goal of computing a point $ \pi $ which is feasible for CP and has better objective value than  $ e$ in that 
\begin{equation}  \label{eqn.bd}
     \frac{c \cdot \pi - z^* }{c \cdot e - z^* } \leq \epsilon \; , 
     \end{equation} 
where $ 0 < \epsilon < 1 $ is user-chosen (i.e., is input to algorithms).  Thus, for the problem of primary interest, CP, the focus is on relative improvement in the objective value.  The reason for focus on relative improvement is due to the following proposition, where the relative error (\ref{eqn.bd}) is shown to be equivalent to an absolute error in approximating the optimal value $ \lambda_{\min,z}^* $ of the CP-equivalent problem (\ref{eqn.bc}).  
In solving a smoothed version of the equivalent problem, the target error needs to be expressed absolutely, in order to make use of the traditional literature on first-order methods. 
\vspace{1mm}

\begin{prop} (\cite[Prop 2.5]{renegar2016efficient}) \label{prop.be}
If $ x \in \mathrm{Affine}_z $ and $ 0 < \epsilon < 1 $,  then 
\begin{align*} 
   & \frac{ c \cdot  \pi(x)  - z^* }{c \cdot e - z^* } \, \leq \,  \epsilon   
    \\ & \qquad   \qquad  \qquad   \textrm{if and only if}    \\ & \qquad  \qquad  \qquad  \qquad   
    \lambda_{\min,z}^* - \lambdamin( x) \, \leq \, \frac{\epsilon }{1 - \epsilon } \, \, \frac{c \cdot e - z }{\, \, \, c \cdot e - z^* } \; .  
\end{align*} 
\end{prop}
 \vspace{1mm}
 
 We conclude development of the basic theory by noting that the theory holds for convex conic optimization problems generally, such as problems of the form
\[ 
  \begin{array}{rl}
  \min_{x \in {\mathcal E}}  & c \cdot x \\
   \textrm{s.t.} & x \in \mathrm{Affine} \\
    & Ax + b \in {\mathcal K}' \end{array}   
\]  
(here, $ A $ is a linear operator from $ {\mathcal E} $ to a Euclidean space $ {\mathcal E}' $, $ b \in {\mathcal E}' $ and $ {\mathcal K}' $ is a proper, closed, convex cone in $ {\mathcal E}' $ with nonempty interior). In the present paper as in \cite{renegar2016efficient}, we focus on problems of the form $ \mathrm{CP} $ because notationally the form is least cumbersome. However, for every result derived in the following sections, an essentially identical result holds for any hyperbolic program regardless of the form, and is even identical in the specific constants.\vspace{1mm}

\subsection{Hyperbolic Programming}  \label{sect.bb}

Here we give a concise introduction to hyperbolic programming, easily accomplished by relying on the reader's grasp of the basics for semidefinite programming. \vspace{1mm}

A hyperbolic program is a convex conic optimization problem for which the underlying cone is a so-called ``hyperbolicity cone.''  The quintessential hyperbolicity cone is $ \Symp  $,  the set of $ n \times n $ positive semidefinite matrices.  The key role in qualifying $ \Symp $ as a hyperbolicity cone is played by the function $ X \mapsto \det(X) $, which if the coefficients of $ X $ are viewed as variables, is a real polynomial of degree $ n $ on the Euclidean space $ \Sym $ ($ n \times n $ symmetric matrices).  The critical property possessed by this polynomial is that there exists $ E \in \mathrm{int}(\Symp) $ -- say, $ E = I $ -- such that for every $ X \in \Sym $, the univariate polynomial $ X \mapsto \det(X - \lambda I) $ has only real roots, the eigenvalues of $ X $. \vspace{1mm}

More generally, a closed convex cone $ {\mathcal K} \subset {\mathcal E}   $ with nonempty interior is said to be a {\em hyperbolicity cone}\footnote{In the literature, the term ``hyperbolicity cone'' generally refers to the interior of $ {\mathcal K}  $, but for us it is expeditious to use the term for the closed cone $ {\mathcal K} $.} if there exists a real  homogeneous polynomial $ p: {\mathcal E} \rightarrow \mathbb{R} $  for which
\begin{enumerate}

\item $ \mathrm{int}({\mathcal K}) $ is a connected component of the set $ \{ x: p(x) \neq 0 \} $ \\
$ \textrm{~} $ \quad (i.e., $ p $ is identically zero on the boundary of $ {\mathcal K} $ \\
$ \textrm{~} $ \qquad \qquad     and nowhere zero in the interior of $ {\mathcal K} $), and

\item for \underline{some} $ e \in \int( {\mathcal K}) $ and all $ x \in {\mathcal E} $ \\
$ \textrm{~} $ \qquad  -- equivalently, by \cite{garding1959inequality}, for \underline{every} $ e \in \int( {\mathcal K}) $ and all  $ x \in {\mathcal E} $ -- \\
$ \textrm{~} $ \quad  the univariate polynomial $ \lambda \mapsto p( x - \lambda e) $ has only real roots.
\end{enumerate}
We say ``$ {\mathcal K} $ has hyperbolic polynomial  $ p $'' (but it should be noted $ p $ is not unique -- for example, $ p^2 $ also is a hyperbolic polynomial for $ {\mathcal K} $).   \vspace{1mm}

It is easy to see that the non-negative orthant, $ \mathbb{R}^n_+ $, is a hyperbolicity cone with polynomial $ p(x) = x_1 \cdots x_n $.  Similarly, it is easy to see the second-order cone $  \{ (x,r) \in \mathbb{R}^n \times \mathbb{R}: \| x \|_2 \leq r \} $ has hyperbolic polynomial $ (x,r) \mapsto r^2 - \| x \|_2^2 $. \vspace{1mm}

In requirement (2), the characterization ``for some'' is useful for verifying hyperbolicity of a given cone, whereas the characterization ``for every'' is critical to easily establishing desirable general properties of hyperbolicity cones.  G\aa rding (\cite{garding1959inequality})  proved the equivalence in the late 1950's, initiating the study of hyperbolicity cones, but in the context of hyperbolic partial differential equations, not optimization.  \vspace{1mm}

Using the characterization ``for every,'' for example, it is simple to verify that if $ {\mathcal K}_1 $ and $ {\mathcal K}_2 $ are hyperbolicity cones whose interiors contain a common point, then $ {\mathcal K}_1 \cap {\mathcal K}_2 $ is a hyperbolicity cone, for which a hyperbolic polynomial is the product $ p = p_1 p_2 $.  Similarly, it is easy to verify that if $ T: {\mathcal E} \rightarrow \bar{{\mathcal E}} $ is a linear transformation whose image intersects the interior of a hyperbolicity cone $ \bar{{\mathcal K}} $ having hyperbolic polynomial $ \bar{p}$, then the pullback $ {\mathcal K} := \{ x \in {\mathcal E}: T(x) \in \bar{{\mathcal K}} \} $ is a hyperbolicity cone, with hyperbolic polynomial $ p(x) = \bar{p}(Tx) $. \vspace{1mm}

Other examples of hyperbolicity cones include the so-called ``derivative cones,'' which form a nested sequence of relaxations of an initial hyperbolicity cone (\cite{renegar2006hyperbolic}). \vspace{1mm}

A {\em hyperbolic program} is a convex conic optimization problem for which the cone is a hyperbolicity cone.  Since we are focusing -- for notational expediency -- on conic optimization problems having the form of CP, for us a hyperbolic program is an optimization problem
\[ \left. \begin{array}{rl}
\min & c \cdot x \\
\textrm{s.t.} & x \in \mathrm{Affine} \\
& x \in {\mathcal K} \end{array}   \right\} \, \, \mathrm{HP} \]
where $ {\mathcal K} $ is a hyperbolicity cone.   \vspace{1mm}

Due to the calculus of hyperbolicity cones (intersections, Cartesian products, pullbacks, etc.), a broad range of optimization problems can be immediately recognized as hyperbolic programs.  This is the flexibility to which we alluded in the introduction when writing that a primary reason for our focus on hyperbolic programming is to make the results -- especially our computational recipes -- be readily applicable in practice (without having to rewrite a model as, say, a semidefinite program). \vspace{1mm}

It was G\"{u}ler \cite{guler1997hyperbolic}  who initiated the study of \vspace{1mm} hyperbolic programming, after realizing that the function $ x \mapsto -\ln p(x) $ is a self-concordant barrier for $ {\mathcal K}  $ (where $ p $ is a hyperbolic polynomial for $ {\mathcal K}  $), and consequently, the general interior-point method theory of Nesterov and Nemirovski \cite{nesterov1994interior} applies.   We recommend \cite{bauschke2001hyperbolic}  as a particularly elegant work examining hyperbolicity cones from the viewpoint of optimization.  \vspace{1mm}   

Letting $ {\mathcal K} $ be a hyperbolicity cone with polynomial $ p: {\mathcal E} \rightarrow \mathbb{R}   $, and fixing $ e \in \mathrm{int}({\mathcal K}) $, for $ x \in {\mathcal E} $ we refer to the roots of $ \lambda \mapsto p(x - \lambda e) $ as the ``eigenvalues of $ x $,'' and denote them $ \lambda_j(x) $ for $ j = 1, \ldots, n $, where $ n := \deg(p) $.  Thus, we list each eigenvalue according to its multiplicity (but we do not assume the eigenvalues are ordered). \vspace{1mm}

The terminology ``eigenvalues'' is a mnemonic device for emphasizing the parallels of the general setting with the special case in which $ {\mathcal K} = \Symp  $ and $ p(X) = \det(X) $. \vspace{1mm}

The eigenvalues of $ x $ depend on the choice of $ e $, but we do not introduce notation making explicit the dependence, because for us $ e $ will be fixed  -- in particular, $ e $ will be the ``distinguished direction'' of the basic theory (\S\ref{sect.ba}).  Observe, then, 
\[ \min_j \lambda_j(x) = \lambda_{\min}(x) \; , \]
where $ \lambda_{\min}(x) $ is as in the basic theory. Also observe, by homogeneity of $ p $,
\begin{equation} \label{eqn.bf} 
 \lambdamax(x) := \max_j \lambda_j (x) = - \lambdamin(-x) \; . 
 \end{equation}

Recall the seminorm $ \| \, \, \|_{\infty} $, defined by (\ref{eqn.bb})  so as to have the property that $ \bar{B}_{\infty}(e,1) $ -- the closed unit ball centered at $ e $ -- is the largest set which both is contained in $ {\mathcal K}  $ and has symmetry point $ e $.  Since for a vector $ v $, 
\[ e + v \in {\mathcal K} \textrm{ and }  e - v \in {\mathcal K} \quad \Leftrightarrow \quad 1 + \lambda_{\min}(v) \geq 0 \textrm{ and } 1 + \lambda_{\min}(-v)  \geq 0 \; , \]
and since $ \lambda_{\min}(-v) = - \lambda_{\max}(v) $, it follows for all $ u \in {\mathcal E} $, 
\begin{equation}  \label{eqn.bg} 
     \| u \|_{\infty} = \max_j | \lambda_j(u)| \; . 
\end{equation}

\section{{\bf  Smoothing}}   \label{sect.c}

The basic theory (\S\ref{sect.ba}) reveals that a convex conic optimization problem is equivalent to a convex optimization problem whose only constraints are linear equations:
\begin{equation}   \label{eqn.ca} 
 \begin{array}{rl}
    \min & c \cdot x \\
    \textrm{s.t.} & x \in \mathrm{Affine} \\
     & x \in {\mathcal K} \end{array} \quad \equiv_e \quad 
     \begin{array}{rl}
       \max & \lambda_{\min}(x) \\
        \textrm{s.t.} & x \in \mathrm{Affine}_z \; , \end{array} \end{equation}
where $ \mathrm{Affine}_z := \{ x \in \mathrm{Affine}: c \cdot x = z \} $ and $ z $ is assumed to satisfy $ c \cdot z < c \cdot e $, with $ e $ being a strictly-feasible point, serving as the distinguished direction. The subscript on the equivalence sign ($ \equiv_e $) is to emphasize that for the equivalence to be useful, $ e $ must be a {\em known} strictly-feasible point. \vspace{1mm}

  In this section it is shown that if $ {\mathcal K} $ is a hyperbolicity cone, then the  function $ x \mapsto \lambda_{\min}(x) $ can be approximated by an explicit, concave, smooth function.  If, for the problem on the right in (\ref{eqn.ca}), the objective function is replaced by the smooth function, we have an optimization problem to which accelerated first-order methods can be applied, so as to approximately solve the problem on the left, the conic problem of interest.  \vspace{1mm}

\subsection{The smooth function and its gradients.}  \label{sect.ca}

To conform to the emphasis on convex (rather than concave) functions in the literature, in this subsection we focus on smoothing the convex function $ x \mapsto \lambdamax(x) $.  Due to the relation (\ref{eqn.bf}), the results easily translate to smoothing the concave function $ x \mapsto \lambdamin(x) $. \vspace{1mm}

Assume $ {\mathcal K} $ is a hyperbolicity cone, with hyperbolic polynomial $ p $ of degree $ n $, and let $ e \in \mathrm{int}( {\mathcal K}) $ be the distinguished direction, giving rise to eigenvalues  $ \lambda_j(x) $ ($ j = 1, \ldots, n $) for every $ x \in {\mathcal E} $. \vspace{1mm}

Our approach is entirely motivated by Nesterov \cite{nesterov2007smoothing}, who focused on the case $ {\mathcal K} = \Symp $. \vspace{1mm}

For $ \mu > 0 $, define the function $ f_{ \mu}: {\mathcal E} \rightarrow \mathbb{R} $ by 
\[ f_{\mu}(x) := \mu \, \ln \sum_j \exp ( \lambda_j(x)/ \mu  )  \; . \]  
It is easily verified for all $ x \in {\mathcal E} $ that
\begin{equation} \label{eqn.cb} 
      \lambda_{\max}(x)  \leq f_\mu(x) \leq \lambda_{\max}(x) + \mu \ln n \; , 
      \end{equation} 
where $ n $ is the degree of the hyperbolic polynomial $ p $ (and hence is the number of eigenvalues, counting multiplicities). \vspace{1mm}

Dual to the seminorm $ \| \, \, \|_{\infty} $ is the extended-real-valued convex function
\[   \| g \|_{\infty}^* := \sup_{\| v \|_{\infty} \leq 1} \lin g, v \rin \;  \]
(this is a norm when $ \| \, \, \|_{\infty} $ is a norm, i.e., when $ {\mathcal K} $ is a regular cone). 
\vspace{1mm} 

\begin{prop}   \label{prop.ca} 
The function $ f_{\mu} $ is convex and analytic.  Moreover, 
\[   \| \grad f_{\mu}(x) - \grad f_{\mu}(y) \|_{\infty}^* \, \leq \, \smfrac{1}{\mu } \, \| x - y \|_{\infty} \quad \textrm{for all $ x,y \in {\mathcal E} $} \; .    \]
\end{prop}
\vspace{3mm}

The proposition is mostly a corollary of Nesterov \cite{nesterov2007smoothing}  and a deep theorem regarding hyperbolicity cones, the Helton-Vinnikov Theorem.  The proof is deferred to Appendix~\ref{app.a}. \vspace{1mm}

For $ f_\mu $ to be relevant to practice, its gradients must be readily computable. We now discuss computation of the gradients. \vspace{1mm}

A key observation is that for many problems, the underlying hyperbolicity cone $ {\mathcal K}  $ is either a Cartesian product of hyperbolicity cones $ {\mathcal K}_i $, or is the intersection  of hyperbolicity cones $ {\mathcal K}_i $. In such cases, computing gradients mostly reduces to consideration of the individual cones $ {\mathcal K}_i $. \vspace{1mm}

\begin{prop} \label{prop.cb}  Assume $ {\mathcal K}_1, \ldots, {\mathcal K}_{\ell} $ are hyperbolicity cones in a Euclidean space $ {\mathcal E} $, with hyperbolic polynomials $ p_1, \ldots, p_{\ell} $. Assume $ e \in \int({\mathcal K}_i) $ for all $ i $, and thus $ e $ can serve as distinguished direction for every $ i $, thereby defining for each $ x \in {\mathcal E} $ a set of eigenvalues $ \{ \lambda_{ij}(x): j = 1, \ldots, n_i \} $ (the roots of $ \lambda \mapsto p_i(x - \lambda e) $), where $ n_i $ is the degree of $ p_i $). Let
\[ f_{i,\mu}(x) := \mu \ln \sum_{j=1}^{n_i} \exp(\lambda_{ij}(x)/\mu ) \; . \]

Let $ {\mathcal K} := \cap_i {\mathcal K}_i $, for which a hyperbolic polynomial is the product $ p = p_1 \cdots p_{\ell} $. For every $ x \in {\mathcal E}$, the corresponding set of eigenvalues -- the roots of $ \lambda \mapsto p(x - \lambda e) $ -- is $ \cup_{i=1}^{\ell} \{ \lambda_{ij}(x): j = 1, \ldots, n_i \} $, and the corresponding smooth function $ f_{\mu} $ satisfies 
\begin{equation}  \label{eqn.zzzz}
   \grad f_{\mu}(x) = \smfrac{1}{\sum_{i=1}^{\ell} \left(  \sum_{j=1}^{n_i} \exp( \lambda_{ij}(x)/\mu) \right) } \,  \sum_{i=1}^{\ell} \left(  \sum_{j=1}^{n_i} \exp( \lambda_{ij}(x)/\mu) \right)   \, \grad f_{i,\mu}(x) \; ,  
   \end{equation} 
a convex combination of the gradients $ \grad f_{i,\mu}(x) $.
\end{prop}
\noindent {\bf Proof:}  The only claim that perhaps is not obvious is the equality (\ref{eqn.zzzz}).  However, the reasoning for the identity becomes simple upon observing the functions $ f_{i, \mu} $ are smooth (Proposition~\ref{prop.ca}), and hence so are the functions  $ x \mapsto \sum_{j=1}^{n_i} \exp( \lambda_{ij}(x)/ \mu ) $.  Consequently,
\begin{align*}
  \grad f_{\mu}(x) &  = \smfrac{\mu}{\sum_{i=1}^{\ell} \sum_{j=1}^{n_i} \exp( \lambda_{ij}(x)/\mu)} \, \sum_{i=1}^{\ell} \grad \left( \sum_{j=1}^{n_i} \exp( \lambda_{ij}(x)/\mu) \right) \\
   & = \smfrac{1}{\sum_{i=1}^{\ell} \sum_{j=1}^{n_i} \exp( \lambda_{ij}(x)/\mu)} \,  \sum_{i=1}^{\ell} \left(  \sum_{j=1}^{n_i} \exp( \lambda_{ij}(x)/\mu) \right)   \, \grad f_{i,\mu}(x) \; , 
   \end{align*}
establishing the identity. \hfill $ \Box $
 \vspace{2mm}

Similarly, if $ {\mathcal K} $ is the Cartesian-product of cones $ {\mathcal K}_i \subseteq {\mathcal E}_i $ with hyperbolic polynomials $ p_i $ ($i = 1, \ldots, \ell $), then $ p(x_1, \ldots, x_m) := p_1(x_1) \cdots p_{\ell}(x_{\ell}) $ is a hyperbolic polynomial for $ {\mathcal K} $.  Moreover, if $ e = (e_1, \ldots, e_{\ell}) \in \int({\mathcal K}) $ is the distinguished direction defining eigenvalues for points $ x = (x_1, \ldots, x_{\ell}) $ in the Euclidean space $ {\mathcal E} = {\mathcal E}_1 \times \cdots \times {\mathcal E}_{\ell} $, then the smooth function $ f_{\mu} $ satisfies
\[  \grad f_{\mu}(x) = (g_1, \ldots, g_{\ell}) \quad \textrm{where }   g_i = \smfrac{\sum_{j=1}^{n_i} \exp(\lambda_{ij}(x_i)/\mu)}{\sum_{k=1}^{\ell} \sum_{j=1}^{n_k} \exp( \lambda_{kj}(x_k)/\mu)} \, \grad f_{i,\mu}(x_i) \; , \]
and where $ \{ \lambda_{ij}(x_i): j = 1, \ldots, n_i \} $ is the set of roots of $ \lambda \mapsto p_i(x_i - \lambda e_i) $. 
\vspace{3mm}

Now we turn to establishing a formula for computing $ \grad f_{\mu}(x) $ once the eigenvalues $ \lambda_j(x) $ are known (to high precision).  In the context of Proposition~\ref{prop.cb}, the identity is pertinent to computing the gradients $ \grad f_{i,\mu}(x) $.
\vspace{1mm}

    The identity is easily motivated in the case that each eigenvalue $ \lambda_j(x) $ is simple (i.e., is a root of multiplicity one for the  \vspace{1mm} polynomial $ \lambda \mapsto p(x - \lambda e) $).

If $ \lambda_j(x) $ is a simple eigenvalue, then  $ y \mapsto \lambda_j(y) $ is analytic in a neighborhood of $ x $ (because simple roots of a univariate polynomial are analytic in small perturbations of the coefficients, and because the coefficients of $ \lambda \mapsto p(x - \lambda e) $ are analytic in $ x $). In particular, when each eigenvalue $ \lambda_j(x) $ is simple, we have 
\begin{equation}   \label{eqn.cc} 
   \grad f_{\mu}(x) = \smfrac{1}{\sum_j \exp(\lambda_j(x)/\mu)} \sum_j \exp(\lambda_j(x)/\mu) \, \,  \grad \lambda_j(x) \; .
   \end{equation}

 On the other hand, differentiating both sides of the equation $ p(y - \lambda_j(y)e) = 0 $ and evaluating at $ x $ leads to
 \begin{equation}  \label{eqn.cd} 
  \grad \lambda_j(x) = \smfrac{1}{p^{(1)}(x - \lambda_j(x)e)} \, \grad p(x - \lambda_j(x)e) \; ,  
  \end{equation}
 where we define for all $ y \in {\mathcal E} $, 
 \[  p^{(1)}(y) :=  \left. \smfrac{d}{dt} p(y + te) \right|_{t=0} \; , \]
a directional derivative. (The denominator in (\ref{eqn.cd})  is nonzero because, by assumption, $ \lambda_j(x) $ is a simple eigenvalue -- a simple root of $ \lambda \mapsto p(x - \lambda e) $ -- and hence $ t = - \lambda_j(x) $  is a simple root of $ t \mapsto p(x + te) $.) \vspace{1mm}

Substituting (\ref{eqn.cd})  into (\ref{eqn.cc})   gives an expression for $ \grad f_{\mu}(x) $ that can be evaluated once the eigenvalues $ \lambda_j(x) $ have been computed (to high precision):
\begin{equation}  \label{eqn.ce} 
 \grad f_{\mu}(x) =  \smfrac{1}{\sum_j \exp(\lambda_j(x)/\mu)} \sum_j \smfrac{\exp(\lambda_j(x)/\mu)}{ p^{(1)}( \,   x - \lambda_j(x)e \, ) } \, \grad p(x - \lambda_j(x)e ) \; ,
 \end{equation}
 a convex combination of the vectors $ \frac{1}{p^{(1)}( \,  x - \lambda_j(x)e \, ) } \, \grad p( \, x - \lambda_j(x)e \, ) $ (for $ j = 1, \ldots, n $, where $ n $ is the degree of $ p $), vectors independent of $ \mu $. \vspace{1mm}
 
 In order to provide an identity similar to (\ref{eqn.ce})  which holds even when eigenvalues are not simple, we introduce the (directional-)derivative polynomials 
\[  p^{(m)}(y)  := \smfrac{d^m}{dt^m} p(y + te)|_{t=0} \; , \quad \textrm{for $ m = 0, 1, \ldots, n $} \; , \]
where $ n $ is the degree of $ p $.  (Like $ p $, these are hyperbolic polynomials for which $ e $ can serve as distinguished direction.  Moreover, using the (inverse) Discrete Fourier Transform to perform interpolation, these polynomials and their gradients can be evaluated efficiently if $ p $ and its gradients can be evaluated efficiently (see \cite{renegar2006hyperbolic}).) \vspace{1mm}

\begin{prop} \label{prop.cc} 
For $ x \in {\mathcal E} $, assume $ \{ \lambda_j(x) \} $ is the set of \underline{distinct} eigenvalues of $ x $, and let $ m_j $ denote the multiplicity of $ \lambda_j(x) $ (that is, the multiplicity of $ \lambda_j(x) $ as a root of $ \lambda \mapsto p(x - \lambda e) $).  Then
\[   \grad f_{\mu}(x) = \smfrac{1}{\sum_j m_j \exp(\lambda_j(x)/\mu)} \sum_j \smfrac{m_j \exp(\lambda_j(x)/\mu)}{p^{(m_j)}( \, x - \lambda_j(x) e \, )} \, \grad p^{(m_j - 1)}( \, x - \lambda_j(x) e \, ) \; , \]
a convex combination of the vectors $ \smfrac{1}{p^{(m_j)}( \, x - \lambda_j(x) e \, )} \, \grad p^{(m_j - 1)}( \, x - \lambda_j(x) e \, ) $ (vectors which are independent of $ \mu $).
\end{prop}
\vspace{2mm}

The proof of Proposition~\ref{prop.cc}  is deferred to Appendix \ref{app.b}. \vspace{1mm}

As a simple example, consider a second-order cone constraint $ \| Ax + b \|_2 \leq \alpha^T x + \beta $, where $ A $ is an $ m \times n $ matrix and $ b \in \mathbb{R}^{m} $, $ \alpha \in \mathbb{R}^n $, $ \beta \in \mathbb{R} $. Assume $ x = \bar{e} $ satisfies the inequality strictly. \vspace{1mm}

Introducing an additional variable $ r $, the constraint can be put into conic form via a rotated Lorentz cone:
\[ 
  \begin{array}{rl}
 r \, \, = 1 & \\
(x,r) \in {\mathcal K} & := \{ (x,r): \| Ax + rb \|_2^2 \leq r( \alpha^T x + \beta r), \, r \geq 0 \textrm{ and } \alpha^T x + \beta r \geq 0  \} \; .
\end{array} 
\] 
A hyperbolic polynomial for $ {\mathcal K} $ is given by 
\begin{align}   
 p(x,r) & := \smfrac{1}{2} \left( \, \| Ax + rb \|_2^2 -  r ( \alpha^T x + \beta r)   \, \right) \nonumber \\
  & \, \,  = \smfrac{1}{2} \left[ \begin{matrix} x^T & r \end{matrix} \right] \left[ \begin{matrix} A^T A & A^T b - \smfrac{1}{2} \alpha \\ b^T A - \smfrac{1}{2} \alpha^T & \| b \|_2^2 - \beta \end{matrix} \right] \left[ \begin{matrix} x \\ r \end{matrix} \right] \; , \label{eqn.cf} 
\end{align}  
and $ e := (\bar{e},1) $ can serve as the distinguished direction. \vspace{1mm}

Then for $ y = \left[ \begin{smallmatrix} x \\ r \end{smallmatrix} \right] $, the eigenvalues are simply the roots $ \lambda^+(y) $ and $ \lambda^-(y) $ for the quadratic polynomial $ \lambda \mapsto \smfrac{1}{2} (y - \lambda e)^T B (y - \lambda e) $, where $ B $ is the matrix in (\ref{eqn.cf}). Moreover, we have the following easily-verified identities: $ \grad p(y) = By $, $ p^{(1)}(y) = e^T By $,  $ \grad p^{(1)}(y) = B e $ and $ p^{(2)}(y) = e^T B e $. \vspace{1mm}

Fixing $ y \in \mathbb{R}^{n+1} $, let $ \lambda^+ := \lambda^+(y) $,  $ y^+ :=  y - \lambda^+ e $, and define $ \lambda^- $, $ y^- $ similarly.  If $  \lambda^+ \neq \lambda^- $, then
\[  \grad f_{\mu}(y) = \smfrac{1}{\exp(\lambda^+/\mu) + \exp( \lambda^-/ \mu) } \left( \smfrac{\exp( \lambda^+/ \mu) }{ e^T B y^+} B y^+  +  \smfrac{\exp( \lambda^-/ \mu) }{ e^T B y^-} B y^- \right) \; , \]
whereas if $ \lambda^+ = \lambda^- $ (i.e., if $ y $ lies in the subspace spanned by $ e $ and the null space of $ B $), then 
\[  \grad f_{\mu}(y) = \smfrac{1}{e^T B e} Be \; , \]
independent of $ \mu $. \vspace{1mm}

More generally, consider a collection of constraints $ \| A_i x + b_i \|_2 \leq \alpha_i^T x + \beta_i $ ($ i = 1, \ldots, \ell $), and assume $ x = \bar{e}  $ strictly satisfies each inequality.  The constraints can be expressed in conic form  as
\[ \begin{array}{rl}
  r \, \, = 1 & \\
(x,r) \in {\mathcal K} &  := \bigcap_i {\mathcal K}_i \; , \end{array} \]
where $ {\mathcal K}_i :=    \{ (x,r): r \| A_ix + rb_i \|_2^2 \leq r( \alpha_i^T x + \beta_i r), \,  r \geq 0 \textrm{ and }  \alpha_i^T x + \beta_i r \geq 0 \} $.  A hyperbolic polynomial for $ {\mathcal K} $ is given by the product $ p = p_1 \cdots p_{\ell} $, where $ p_i(y) = \smfrac{1}{2} y^T B_i y $ and $ B_i $ is the matrix in (\ref{eqn.cf})  but with entries subscripted.  Again $ e = (\bar{e},1) $ can serve as distinguished direction. \vspace{1mm}

Now $ f_{\mu} $ is a smoothing of the function $ y \mapsto \max \{ \lambda_i^+(y), \lambda_i^-(y): i = 1, \ldots, \ell \} $.  To compute $ \grad f_{\mu}(y) $, one can first perform the simple computations above to compute $ \grad f_{i,\mu}(y) $ for each $ i $, then compute the convex combination of these gradients using the weights prescribed by Proposition~\ref{prop.cb}.   

\subsection{The setup for applying accelerated methods.}  \label{sect.cb}

For $ \mu > 0 $, let $ \hat{f}_\mu: {\mathcal E} \rightarrow \mathbb{R} $ be the smooth \underline{concave} function
\[   \hat{f}_\mu(x) := - f_\mu(- x) \, \, = - \mu \ln \sum_j \exp( - \lambda_j(x)/\mu) \; , \]
an approximation to the function $ x \mapsto \lambdamin(x) $ in satisfying
\begin{equation}  \label{eqn.cg} 
       \lambdamin (x) - \mu \ln n \leq \hat{f}_\mu(x)  \leq \lambdamin(x) \; . 
        \end{equation} 

We have three optimization problems, 
\begin{equation}  \label{eqn.ch} 
   \textrm{HP}  \left\{ \, \begin{array}{rl}
          \min & c \cdot x \\
            \textrm{s.t.} & x \in \mathrm{Affine}  \\
               & x \in {\mathcal K} \end{array} \right.  \, \, \equiv_e \, \, \begin{array}{rl}
               \max & \lambdamin(x) \\
                 \textrm{s.t.} & x \in \mathrm{Affine}_z \end{array} \, \,  \approx_e \, \,  \begin{array}{rl}
   \max & \hat{f}_{\mu}(x) \\
   \textrm{s.t.} & x \in \mathrm{Affine}_z \; , \end{array}   \end{equation}
the first two being equivalent (assuming $ z < c \cdot e $ and $ e $ is a known strictly-feasible point), and the third being a smoothing of the second.  To approximately solve HP, we apply an accelerated method to the third problem, with $ \mu $ chosen appropriately. \vspace{1mm}

In applying an accelerated method, we adopt a Euclidean perspective, endowing $ {\mathcal E} $ with an inner product $ \langle \; , \; \rangle $ to be used for computation, so as to make it easy to retain feasibility for the smoothed problem, simply by orthogonally projecting gradients onto $ {\mathcal L} $, the subspace of which $ \mathrm{Affine}_z $ is a translate.  Let $ P_{{\mathcal L}} $ denote the projection operator.  Of course $ {\mathcal L} $ -- and hence $ P_{{\mathcal L}} $ -- is independent of $ z $. \vspace{1mm}

The inner product $ \langle \; , \; \rangle $ -- the ``computational inner product'' -- can differ from the inner product used in expressing HP, i.e., the inner product $ (u,v) \mapsto u \cdot v $. We make the distinction to highlight that the choice of computational inner product can affect convergence rates, which is not surprising given that supgradients and orthogonal projections are dependent on the inner product.  Precisely how the computational inner product affects the rates is discussed immediately after Corollary 4.3, where it also is noted that except for problems with special structure, choosing the ``best'' inner product is an extremely difficult task. \vspace{1mm}

Let $ \| \, \, \| $ denote the norm corresponding to the computational \vspace{1mm} inner product, $ \| v \| = \lin v, v \rin^{1/2} $.  

The computational inner product provides an identification of $ {\mathcal E} $ with its dual space.  Henceforth we use $ \grad \hat{f}_{\mu}(x)  $ to denote the gradient viewed as an element of $ {\mathcal E} $.\footnote{That is, $ g = \grad \hat{f}_{\mu}(x) $ is the vector in $ {\mathcal E} $  satisfying $ 0 = \lim_{\| \Delta x \| \downarrow 0} \smfrac{\hat{f}_{\mu}(x + \Delta x) - \hat{f}_{\mu}(x) - \lin g, \Delta x \rin}{ \| \Delta x \|   } $.}    
    
In this setup, the relevant Lipschitz constant for accelerated methods is a value $ L_{\mu,z} $ satisfying 
\[  \| P_{{\mathcal L}} \grad \hat{f}_{\mu}(x) \, - \, P_{{\mathcal L}} \grad \hat{f}_{\mu}(y) \| \, \leq \, L_{\mu,z} \, \| x - y \| \quad \textrm{for all $ x, y \in \mathrm{Affine}_z $} \; .  \]
Generally, the best such constant depends on $ z $, but the dependence is mild. Considering the special case $ \mu = 1 $, we are thus led to define a Lipschitz constant applying to all translates of $ {\mathcal L} $: 
\begin{equation}  \label{eqn.ci}  
   L_1 = \sup_{u \in {\mathcal E}} \, \, \sup \{ \smfrac{\| P_{{\mathcal L}} \grad \hat{f}_1 (u + v) - P_{{\mathcal L}} \grad \hat{f}_1(u + w) \|  }{\| v - w \| } : v, w \in {\mathcal L} \textrm{ and }   v \neq w \} \; . 
   \end{equation}
Depending on the subspace $ {\mathcal L} $, the value $ L_1 $ can be far smaller than the Lipschitz constant for the unrestricted gradient map, that is, for the case $ {\mathcal L} = {\mathcal E} $.   
\vspace{1mm}

Observe that  $ \hat{f}_{\mu}(x) = \mu \hat{f}_1( \smfrac{1}{\mu } x) $ (due to positivity of $ \mu $, the definition of $ \hat{f}_{\mu} $, and the homogeneity $ \lambda_j(\smfrac{1}{\mu } x) = \lambda_j(x)/\mu $).  Consequently, if $ \hat{f}_{\mu} $ is substituted for $ \hat{f}_1 $ in (\ref{eqn.ci}), the supremum is equal to $ L_1/ \mu  $. \vspace{1mm}

The value $ L_1 $ both reflects properties of the specific subspace $ {\mathcal L} $ and immediately gives the Lipschitz constant $ L_1/\mu  $ for the gradient map $ x \mapsto P_{{\mathcal L}} \grad \hat{f}_{\mu}(x) $ restricted to any of the affine spaces $ \mathrm{Affine}_z $ (thus eliminating the need for a Lipschitz constant $ L_{\mu, z} $ subscripted with parameters $ \mu $ and $ z $).   

The following corollary to Proposition~\ref{prop.ca}  provides geometric insight into how the constant $ L_1 $ reflects the particular subspace $ {\mathcal L} $ in relation to the hyperbolicity cone $ {\mathcal K} $,  as well as reflects the choices for the distinguished direction $ e $ and the computational inner product $ \langle \; , \; \rangle $ (for which, we recall, the norm is written $ \| \, \, \| $).  \vspace{1mm}

Define 
\begin{equation}  \label{eqn.cj} 
   r_e  = \sup \{ r: B(e,r) \cap \mathrm{Affine}_{c \cdot e} \subseteq {\mathcal K} \} \; ,  
   \end{equation} 
where $ B(e,r) $ is the $ \| \, \, \| $-ball of radius $ r $ centered at $ e $. 

\begin{cor} \label{cor.cd} \quad  $  L_1 \leq 1/r_e^2  \; . $   
\end{cor}
\noindent {\bf Proof:}  As $ \bar{B}_{\infty}(e,1) $ is the largest set both contained in $ {\mathcal K} $ and having symmetry point $ e $, the intersection $ \bar{B}_{\infty}(e,1) \cap \mathrm{Affine}_{c \cdot e} $ is the largest set both contained in $ {\mathcal K} \cap \mathrm{Affine}_{ c \cdot e} $ and having symmetry point $ e $.  Thus, the Euclidean norm satisfies
\[   \bar{B}(e, r_e ) \cap \mathrm{Affine}_z \, \subseteq \, \bar{B}_{\infty}(e,1) \cap \mathrm{Affine}_z \; , \]
from which follows
\[   
   \| \bar{v}  \|_{\infty} \leq \| \bar{v}  \|/ r_e  \quad \textrm{for all $ \bar{v}  \in {\mathcal L} $} \; . 
\] 
Consequently, for $ u \in {\mathcal E} $ and $ v,w \in {\mathcal L} $,  
\begin{align*}
& \| P_{{\mathcal L}} \grad \hat{f_{1}}(u+v) \, - \, P_{{\mathcal L}} \grad \hat{f_{1}}(u+w) \| \\
& = \max \{ \lin \bar{v} ,  \grad \hat{f}_1(u+v) - \grad \hat{f}_1(u+w) \rin: \bar{v}  \in {\mathcal L} \textrm{ and } \| \bar{v}  \| \leq 1 \} \\
& \leq \max \{  \lin \bar{v} ,  \grad \hat{f}_1(u+v) - \grad \hat{f}_1(u+v) \rin: \bar{v} \in {\mathcal E} \textrm{ and }      \| \bar{v}  \|_{\infty} \leq 1/ r_e   \} \\   
& =   \|   \grad \hat{f}_1(u+v) - \grad \hat{f}_1(u+w) \|_{\infty}^*/ r_e   \\
& \leq \| v - w \|_{\infty}/r_e   \quad \textrm{(by Proposition~\ref{prop.ca})} \\
& \leq \| v - w \|/r_e^2  \; ,  
\end{align*}
completing the proof (due to the definition of $ L_1 $). \hfill $ \Box $
 \vspace{3mm}
 
We close the section with a lemma that ties accuracy in solving the smoothed problem to the accuracy in solving HP, the problem of primary interest. \vspace{1mm}

Recall $ z^* $ denotes the optimal value of HP, the optimization problem of primary interest -- the first problem in (\ref{eqn.ch}).  Recall $ \lambda_{\min,z}^* $ is the optimal value of the HP-equivalent problem -- the second problem in (\ref{eqn.ch})  -- and $ x_z^* $ is any optimal solution.  Let $ f_{\mu,z}^* $ denote the optimal value of the smoothed problem -- the third problem -- and let $ x_{\mu, z}^* $ be an optimal solution.

\begin{lemma}  \label{lem.ce} If $ x \in \mathrm{Affine}_z $ where $ z < c \cdot e $, then the projection $ \pi = \pi(x) $ satisfies
\[  \frac{c \cdot \pi - z^*}{c \cdot e - z^*} \leq \left( \hat{f}_{\mu,z}^*  - \hat{f_{\mu}}(x) + \mu \, \ln n \right) \, \frac{c \cdot e - z^*}{c \cdot e - z} \; . \]
\end{lemma}
\noindent {\bf Proof:}  Define $ \epsilon  := \frac{c \cdot \pi - z^*}{ c \cdot e - z^*} $.  According to Proposition~\ref{prop.be}, 
\[  \lambda_{\min,z}^* - \lambda_{\min}(x) = \smfrac{\epsilon}{1 - \epsilon} \, \smfrac{c \cdot e - z}{c \cdot e - z^*} \, \geq \, \epsilon \, \smfrac{c \cdot e - z}{c \cdot e - z^*} \; . \]
Thus, to prove the lemma it suffices to show
\[  \lambda_{\min,z}^* - \lambda_{\min}(x) \leq  \hat{f}_{\mu,z}^*  - \hat{f_{\mu}}(x) + \mu \, \ln n \; . \]
However, using (\ref{eqn.cg}),
\begin{align*}
 \lambda_{\min,z}^* - \lambda_{\min}(x) & = 
 \big( \lambda_{\min,z}^* - \hat{f_{\mu}}(x_z^*) \big) + \big( \hat{f_{\mu}}(x_z^*) - \hat{f}_{\mu,z}^*  \big) \\ & \qquad  + \big( \hat{f}_{\mu,z}^*   - \hat{f_{\mu}}(x) \big) + \big( \hat{f_{\mu}}(x) - \lambda_{\min}(x) \big) \\
& \leq \mu \, \ln n + 0 + \big( \hat{f}_{\mu,z}^*   - \hat{f_{\mu}}(x) \big) + 0 \; , 
 \end{align*}
 completing the proof. \hfill $ \Box $
 \vspace{7mm}

\section{{\bf  Applying Accelerated First-Order Methods}} \label{sect.d}

Virtually any (accelerated) first-order method can be applied to solving 
\begin{equation}  \label{eqn.da} 
  \begin{array}{rl}
   \max & \hatfmu(x) \\
     \textrm{s.t.} & x \in \mathrm{Affine}_z \; ,   \end{array}  
     \end{equation} 
as feasibility can be retained using only orthogonal projection of gradients onto a subspace, the same subspace at every iteration.  We develop results for an accelerated method due to Nesterov \cite{nesterov2013gradient}.  Among the accelerated methods not requiring (an upper bound on) the Lipschitz constant as input, this algorithm streamlines most nicely when specialized to the setting where the only projections are onto a subspace.  However, similar complexity results could be obtained if instead we employed the universal method of Lan \cite{lan2015bundle}, or of Nesterov \cite{nesterov2014universal}, and would be obtained if we employed any of numerous accelerated methods that {\em  do} require a Lipschitz constant as input (an unrealistic requirement in our general setting). \vspace{1mm}

 Following is the accelerated method \cite[(4.9)]{nesterov2013gradient} specialized to our setting (see Appendix \ref{app.c}  for discussion).    To simplify specification of the algorithm, the strictly-feasible point $ e $ is viewed as fixed, not input. 

\noindent 
\hrulefill
 
\noindent {\bf AGM} (Accelerated Gradient Method) \vspace{1mm}

\noindent  
(0) Input:  $ L_1^{\bullet} > 0 $, an estimate of the value $ L_1 $ defined by (\ref{eqn.ci}), \\
$ \textrm{~} $ \qquad \qquad   $ \mu > 0 $, and \\
$ \textrm{~} $ \qquad \qquad     $ x_0  \in \mathrm{Affine} $ satisfying $ z := c \cdot x_0  < c \cdot e $. \\
$ \textrm{~} $ \quad Initialize:   $ L = L_1^{\bullet}/\mu   $, \, $ v = x_0 $, \, $ A = 0 $ \,
and \,
 $ k = 0 \; .   $    \\
(1) Compute $ a = \frac{1 + \sqrt{1 + 2 AL}}{L} $, \,
 $ y = \frac{A}{A + a} x_k + \frac{a}{A + a} v $ \,
  and \,
   $ x = y + \frac{1}{L} P_{{\mathcal L}} \grad \hat{f_{\mu}}(y) \; . $    \\  
(2) If $ \lin \grad \hat{f_{\mu}}(x), y - x \rin > - \frac{1}{L} \,  \| P_{ {\mathcal L}} \grad \hat{f_{\mu}}(x) \|^2 $, let $ L \leftarrow 2L $ and return to Step 1. \\
(3) Define $ x_{k+1} = x  $ and let \, $ v \leftarrow v + a \, P_{{\mathcal L}} \grad \hat{f_{\mu}}(x)  $,  \, $ L \leftarrow L/2 $, \, $ A \leftarrow A + a \; , $     \\
$ \textrm{~} $ \qquad \qquad  \qquad  \qquad  \qquad \qquad  \qquad  \qquad  \qquad  \qquad  \qquad  \,   $ k \leftarrow k+1 $, then   \vspace{-1.5mm} go to Step 1.

\noindent 
\hrulefill
\vspace{2mm}

An important aspect of AGM is that in moving from iterate $ x_k $ to $ x_{k+1} $, numerous gradient steps $ x = y + \smfrac{1}{L} P_{\mathcal L} \grad \hat{f}_{\mu}(y) $   might have to be tried before arriving at a step which is not rejected in (2).  We are thus careful to distinguish between the number of iterates $ x_k $ in a sequence and the number of gradients computed when constructing the sequence. \vspace{1mm}

It is easily seen that all iterates $ x_k $ lie in $ \mathrm{Affine}_z $ for $ z = c \cdot x_0 $. \vspace{1mm} 

Recall $ \hat{f}_{\mu,z}^* $ denotes the  optimal value of (\ref{eqn.da}), and $ x_{\mu,z}^* $ is any optimal solution.
\vspace{1mm}

\begin{thm}  {\bf (Nesterov)}  \label{thm.da} 
Let $ z = c \cdot x_0 $ and $  
  k^{\bullet}  := \max \{ 0, \lceil \,  \log_2 (L_1^{\bullet}/L_1) \,  \rceil \} \; . $ 
  Then 
\[  k > k^{\bullet} \quad \Rightarrow \quad   \hat{f}_{\mu,z}^* - \hat{f_{\mu}}(x_k) \, \leq \, \frac{2 L_1}{\mu} \,     \left(  \frac{ \| x_0 - x_{\mu, z  }^* \|}{ k-k^{\bullet}} \right)^2  \; . \]
Moreover, for $ k > k^{\bullet} $, the total number of gradients computed before $ x_k $ is determined does not exceed 
\begin{equation} \label{eqn.db} 
  4 (k - k^{\bullet}) + 2  \left| \log_2(  L_1^{\bullet}/L_1   ) \right| + 6 \; . 
  \end{equation}   
\end{thm}
\vspace{2mm}

Results from \cite{nesterov2013gradient}  can readily be assembled into a proof of Theorem~\ref{thm.da}, as is demonstrated in Appendix \ref{app.c}. 
 \vspace{1mm}
 
For recasting Theorem~\ref{thm.da}  into a form meaningful to HP (the optimization problem of primary interest), we follow the approach introduced in \cite{renegar2016efficient}, to which the present paper is a sequel. \vspace{1mm}

 The {\em level sets}  for HP are the sets
\[  \mathrm{Level}_{z} =   \mathrm{Affine}_{z} \cap {\mathcal K}   \; ,   \]
that is, the largest feasible sets for HP on which the objective function is constant\footnote{There is possibility of confusion here, as in the optimization literature, the terminology ``level set'' is often used for the portion of the feasible region on which the (convex) objective function does not exceed a specified value rather than -- as for us -- exactly equals the value.  Our choice of terminology is consistent with the general mathematical literature, where the region on which a function does not exceed a specified value is referred to as a sublevel set, not a level set.}.  If $ z < z^* $, then $ \mathrm{Level}_{z} = \emptyset \; . $    \vspace{1mm}

If some level set is unbounded, then either HP has unbounded optimal value or can be made to have unbounded optimal value with an arbitrarily small perturbation of $ c $.  In the same vein, a level set with large diameter is an indication that perturbations can result in relatively large changes in optimal value.   \vspace{1mm}  
     
For scalars $ z  $, let 
\[    \diamval  :=  \sup  \{  \| x - y \|: x,y \in \mathrm{Level}_{z}   \} \; ,    \]
the diameter of $ \mathrm{Level}_{z} $.   (If $ \mathrm{Level}_z = \emptyset $, let $ \diamval := - \infty  $.) \vspace{1mm} 

Recall $ \lambda_{\min,z}^* $ denotes the optimal value of the HP-equivalent problem (the problem to which (\ref{eqn.da})  is a smoothing). \vspace{1mm} 
 
\begin{cor} \label{cor.db} 
Assume input $ x_0 $ is strictly feasible for $ \mathrm{HP} $, and let $ z := c \cdot x_0 $. Assume $ \epsilon'  $ is a scalar satisfying $ 0 < \epsilon' \leq 2 \lambda_{\min}(x_0) $. If input $ \mu $ is chosen as $ \mu =  \smfrac{1}{2} \epsilon'/ \ln(n) $, then 
\[  k > k^{\bullet} +  \frac{2 \, \mathrm{diam}_z  \sqrt{2 L_1  \ln n}}{\epsilon'} \quad \Rightarrow \quad \lambda_{\min,z}^* - \lambda_{\min}(x_k) \leq \epsilon' \; , \]
where $ k^{\bullet} :=  \max \{ 0, \lceil  \, \log_2 (L_1^{\bullet}/L_1) \, \rceil \} \; . $ \vspace{1mm}

Moreover, for $ k > k^{\bullet} $, the total number of gradients computed before $ x_k $ is determined does not exceed \, $ 4 (k - k^{\bullet}) + 2 \left| \log_2(  L_1^{\bullet}/L_1   ) \right| + 6 \; . $ 
\end{cor}
\noindent {\bf Proof:}  Fix $ \mu = \smfrac{1}{2} \epsilon'/ \ln(n) $. It is readily seen from (\ref{eqn.cg})  that for $ x \in \mathrm{Affine}_z $, 
\[  \hat{f}_{\mu,z}^* - \hat{f}_{\mu}(x) \leq \smfrac{1}{2} \epsilon' \quad \Rightarrow \quad \lambda_{\min,z}^* - \lambda_{\min}(x) \leq \epsilon' \; . \]
Thus, Theorem~\ref{thm.da}  is easily seen to imply
\[  k > k^{\bullet} + \frac{2 \| x_0 - x_{\mu,z}^* \| \sqrt{ 2 L_1 \ln n}}{ \epsilon'} \quad \Rightarrow \quad \lambda_{\min,z}^* - \lambda_{\min}(x_k) \leq \epsilon' \; . \]
To complete the proof, it suffices to show $ x_{\mu,z}^* $ is feasible for HP, as then the inequality \\
 $    \| x_0 - x_{\mu,z}^* \| \leq \mathrm{diam}_z $ 
follows from the assumed feasibility of $ x_0 $. \vspace{1mm}

However, since $ 0 < \epsilon' \leq 2 \lambda_{\min}(x_0) $, we have 
\[  \mu = \smfrac{1}{2} \epsilon'/ \ln(n) \leq \lambda_{\min}(x_0)/ \ln(n) \; , \]
and thus (\ref{eqn.cg})  implies $ \hat{f}_{\mu}(x_0) \geq 0 $.  Hence, $ \hat{f}_{\mu}(x_{\mu,z}^*) \geq 0 $, which by (\ref{eqn.cg})  implies $ \lambda_{\min}(x_{\mu,z}^*) \geq 0 $, that is, implies $  x_{\mu,z}^* $ is indeed feasible for HP. \hfill $ \Box $
 \vspace{2mm}

For the iterates $ x_k $  of the accelerated method, define $ \pi_k = \pi(x_k) $.  The main algorithmic goal is achieved when $ \pi_k $ is computed for which $ \frac{c \cdot \pi_k - z^*}{c \cdot e - z^*} \leq \epsilon $.  Proposition~\ref{prop.be}  shows the goal is accomplished precisely when iterate $ x_k $ satisfies 
 \begin{equation}  \label{eqn.dc} 
            \lambda_{\min,z}^* - \lambda_{\min}(x_k) \, \leq \, \frac{\epsilon}{1 - \epsilon} \, \,   \frac{c \cdot e - z}{c \cdot e - z^*} \; . 
  \end{equation}
It is thus tempting to try making use of Corollary~\ref{cor.db}  with $ \epsilon' $ chosen to be the value on the right of (\ref{eqn.dc}).  However, this is an undesirable approach for multiple reasons, one being that the quantity on the right of (\ref{eqn.dc})  depends on the optimal value $ z^* $ (which we want to avoid assuming known), and another being that the resulting bound on the total number of gradients computed would be of the form
\begin{equation}  \label{eqn.dd} 
 {\mathcal O} \left( \, \frac{ \mathrm{diam}_z \sqrt{L_1 \ln n}}{ \epsilon} \, \cdot \, \frac{c \cdot e - z^*}{ c \cdot e - z} \, + \, \left|  \log    \frac{L_1^{\bullet}}{L_1}   \right| \, \right) \; , 
 \end{equation}
 If $ z $ -- the objective value for the known feasible point $ x_0 $ -- is only slightly smaller than $ c \cdot e $, then the ratio $ \frac{c \cdot e - z^*}{c \cdot e - z} $ can be large compared to $ 1/ \epsilon $, in which case the bound (\ref{eqn.dd})  is less appealing.  \vspace{1mm}

To arrive at an algorithm which does not require $ z^* $ to be known, and for which the disconcerting product $ \frac{1 }{\epsilon} \, \frac{c \cdot e - z^*}{c \cdot e - z} $ does not appear as it does in (\ref{eqn.dd}),  we must go beyond simply applying an accelerated method to the smoothed problem (\ref{eqn.da}).  Towards presenting a strategy which accomplishes these goals, we first record a result highlighting a case where the ratio $ \frac{c \cdot e - z^*}{c \cdot e - z} $ is small, leading to a bound of the form (\ref{eqn.dd})  which is more appealing. 
\vspace{1mm}

\begin{cor}  \label{cor.dc} 
Assume the input $ x_0 $ is strictly feasible for $ \mathrm{HP} $ and satisfies
\begin{equation}  \label{eqn.de} 
     \lambda_{\min}(x_0) \geq 1/6 \quad \textrm{and} \quad \frac{c \cdot e - z^*}{c \cdot e - z} \, \leq \,  3 \; , 
     \end{equation}
where $ z := c \cdot x_0 < c \cdot e $.  Let $ 0 < \epsilon < 1 $. 

If input $ \mu $ is chosen as $ \mu = \smfrac{1}{6} \epsilon/\ln(n) $, then for the projected iterates $ \pi_k := \pi(x_k) $, 
\begin{equation}  \label{eqn.df} 
   k > k^{\bullet} +  \frac{ 6 \, \mathrm{diam}_z  \sqrt{2 L_1 \ln n}}{\epsilon } \quad \Rightarrow \quad \frac{c \cdot \pi_k - z^*}{c \cdot e - z^*} \, \leq \, \epsilon \; ,    
   \end{equation} 
where \, $ k^{\bullet} := \max \{ 0, \lceil \,   \log_2(L_1^{\bullet}/L_1) \, \rceil  \} \; .  $ \vspace{1mm}

Moreover, for $ k > k^{\bullet} $, the total number of gradients computed before $ x_k $ is determined does not exceed \, $ 4 (k - k^{\bullet}) + 2 \left|  \log_2(  L_1^{\bullet}/L_1  ) \right| + 6 \; . $ 
\end{cor}  
\noindent {\bf Proof:}  
Let $ \epsilon' := \epsilon/3 \leq 1/3 $.  Then, by assumption, $ \epsilon' \leq 2 \lambda_{\min}(x_0) $, and hence Corollary~\ref{cor.db} can be applied with $ \mu = \smfrac{1}{2} \epsilon'/\ln(n) = \smfrac{1}{6} \epsilon/ \ln(n) $,
giving 
\begin{align*}
    k > k^{\bullet} +  \frac{ 6 \, \mathrm{diam}_z \sqrt{2 L_1 \ln n}}{\epsilon } & \, \Rightarrow \,  \lambda_{\min,z}^* - \lambda_{\min}(x_k) \leq \epsilon/3 \\
   & \, \Rightarrow \,  
   \lambda_{\min,z}^* - \lambda_{\min}(x_k) \, \leq \frac{\epsilon }{1 - \epsilon} \, \, \frac{c \cdot e - z}{c \cdot e - z^* } \; ,  
    \end{align*}
where the last implication is due to the rightmost inequality assumed in (\ref{eqn.de}).  Invoking Proposition~\ref{prop.be} completes the proof.   \hfill $ \Box $
\vspace{2mm}

The threshold on $ k $ given by (\ref{eqn.df})  can be given geometric meaning by recalling the value $ r_e  $ defined by (\ref{eqn.cj}), a value reflecting the extent to which the computational inner product conforms to the geometry of the set obtained by intersecting the cone $ {\mathcal K}  $ with $ \mathrm{Affine}_{c \cdot e} $.    Substituting $ L_1 \leq 1/r_e^2 $ (by Corollary~\ref{cor.cd})  gives
\[ 
 k > k^{\bullet} +  \frac{6 \, \sqrt{\ln n}}{ \epsilon }  \, \,   \frac{ \mathrm{diam}_z}{ r_e}  \quad \Rightarrow \quad \frac{c \cdot \pi_k - z^*}{c \cdot e - z^*} \, \leq \, \epsilon \; . 
\] 
Here, the strongest result is obtained, of course, if the computational inner product is an inner product which minimizes the ratio $ \mathrm{diam}_z/ r_e $. (Except for problems with special structure, however, determining an inner product that nearly minimizes the ratio is a difficult task.) \vspace{1mm}
  
Corollary~\ref{cor.dc}  is limited in its scope, as the second assumption in (\ref{eqn.de})  is strong. Nonetheless, to obtain an efficient algorithm, we apply the accelerated method sequentially in such a manner that on the final application, the assumptions of Corollary~\ref{cor.dc}  are satisfied.  \vspace{1mm}

Our main algorithm, given below, has both outer iterations and inner iterations.  At the beginning of an outer iteration, two different inputs to AGM are specified, on which the method is run separately, in parallel. To simplify the presentation, it is useful to introduce notation which reflects the relevant choices of input to AGM  -- here, the approximation $ L_1^{\bullet} $ of the Lipschitz constant $ L_1 $ should be viewed as fixed (along with $ e $, the distinguished direction):  
\begin{center}  $ \textrm{AGM}( \mu, v) $ \, -- \,  AGM run with input $ e $, $ L_1^{\bullet} $, $ \mu $ and $ x_0 = v \; . $ 
\end{center}
\vspace{1mm}
Likewise, to simplify specification of the main algorithm, $ L_1^{\bullet} $ (and $ e $) are viewed as fixed, not as input.

\noindent 
\hrulefill
 
\noindent {\bf MainAlgo} \vspace{1mm}

\noindent  
Input: $ 0 < \epsilon < 1  $ (the accuracy desired in solving HP), and \\
$ \textrm{~} $ \qquad \quad     $ v_1     \in \mathrm{Affine} $ satisfying $ c \cdot v_1    < c \cdot e $ and\footnote{For any vector $ v \in \mathrm{Affine} $ satisfying $ c \cdot v < c \cdot e $, the vector $ v_1 = e + \smfrac{3}{4} \, \smfrac{1}{1 - \lambda_{\min}(v) }   ( v - e) $ is appropriate as input.} $ \lambda_{\min}(v_1) = 1/4 $.  
\vspace{1mm}

\noindent
Fix: $ \mu_1 = 1/(12 \ln n) $, and \\
$ \textrm{~} $ \quad   \, $ \mu_2 = \epsilon/(6 \ln n) $.
\vspace{1mm}

\noindent 
Initiate:  $  \ell = 1 $ (counter for outer iterations).
\vspace{2mm}

\noindent  
Procedure: \vspace{1mm}

\noindent 
{\bf  (1)} Start, in parallel, $ \textrm{AGM}_1 :=  \textrm{AGM}(\mu_1, v_{\ell}) $ and $ \textrm{AGM}_2 :=  \textrm{AGM}(\mu_2, v_{\ell}) $, letting the implementations run indefinitely  unless an iterate $ w_j $ for $ \textrm{AGM}_1 $  is obtained for which $ \lambda_{\min}(w_j) \geq 1/2 $, in which case go to Step 2.
\vspace{1mm}

\noindent  
{\bf  (2)} For the iterate $ w_j $ satisfying $ \lambda_{\min}(w_j) \geq 1/2 $, let $ v_{\ell + 1} := e + \smfrac{3}{4} \, \smfrac{1}{1 - \lambda_{\min}(w_j)} ( w_j - e ) $, an ``outer iterate.''  Let $ \ell \leftarrow \ell + 1 $, then go to Step 1, restarting $ \mathrm{AGM}_1 $ and $ \mathrm{AGM}_2 $ with the new $ v_{\ell} $  (but leaving $ \mu_1 $ and $ \mu_2 $ unchanged).   
   
\noindent 
\hrulefill
\vspace{2mm}

In practice, a smarter way to proceed might be to  update the estimate $ L_1^{\bullet} $ based on the last values for $ L $ just before $ \textrm{AGM}_1 $ and $ \mathrm{AGM}_2 $ are restarted.    \vspace{1mm}
 
Denote the sequence of iterates for $ \mathrm{AGM}_1 $ by $ w_1, w_2, \ldots $, where the indexing is {\em not}  restarted when $ \mathrm{AGM}_1 $ is restarted.  (Thus, if $ w_j $ satisfies the property which causes Step 2 to be invoked and the two implementations to be restarted, the next iterate produced by $ \mathrm{AGM}_1 $ is $ w_{j+1} $.)  Likewise, denote the sequence of iterates for $ \mathrm{AGM}_2 $ by $ x_1, x_2, \ldots $, where the indexing is {\em  not} restarted when $ \mathrm{AGM}_2 $ is restarted.    \vspace{1mm}

To establish a bound on the total number of gradients computed, we assume $ \mathrm{AGM}_1 $ and $ \mathrm{AGM}_2 $   compute gradients simultaneously, so that at any given time, both implementations have computed the same number of gradients. \vspace{1mm}

Let $ z_1 := c \cdot v_1 $, where $ v_1 $ is input to MainAlgo, and define
\[  \mathrm{Diam}_{z_1} = \max \{ \mathrm{diam}_z: z \leq z_1 \} \; . \]

\begin{thm} \label{thm.dd}  
$ \mathrm{MainAlgo} $ produces an iterate $ x_k $ for which $ \pi_k := \pi( x_k) $ satisfies $ \frac{c \cdot \pi_k - z*}{c \cdot e - z^*} \leq \epsilon $, and does so within computing a total number of gradients  not exceeding
\begin{align*}
&  {\mathcal O}   \left( \, \,  \frac{\mathrm{Diam}_{z_1} \sqrt{ L_1 \ln n}}{\epsilon} \, \right. \\
& \qquad  \quad  \left. + \, \left( 1 + \log \,  \frac{c \cdot e - z^*}{c \cdot e - z_1}  \right) \, \cdot \, \left( \, 1 \, + \,  \mathrm{Diam}_{z_1} \sqrt{ L_1 \ln n}        \, + \,  \left|  \log   \frac{L_1^{\bullet}}{L_1}  \right|  \,  \right) \,   \right) \; . 
\end{align*}                      
\end{thm}   
\vspace{2mm}

\noindent {\em Remark 1:} Unlike the bound (\ref{eqn.dd}), here the ratio $ \smfrac{c \cdot e - z^*}{c \cdot e - z_1} $ has been entirely separated from terms involving $ \epsilon $.  Moreover, now only the logarithm of the ratio appears. \vspace{2mm}

\noindent {\em Remark 2:} Similar to the observation immediately following Corollary~\ref{cor.dc}, if $ 1/r_e^2 $ is substituted for $ L_1 $ in $ \sqrt{L_1 \ln n} $ (Corollary~\ref{cor.cd}), the result is strongest when the computational inner product is chosen to minimize the ratio $ \textrm{Diam}_{z_1}/r_e $. \vspace{2mm}

\noindent {\bf Proof of Theorem~\ref{thm.dd}:}  As noted just before Theorem~\ref{thm.da}, all iterates of AGM have the same objective value as the initial iterate.  Consequently, for sequences of iterates $ w_j $ and $ x_k $ generated by MainAlgo, all iterates have the same objective value so long as Step 2 is not invoked while the sequences are being generated. Let us examine the change in objective value when Step 2 is invoked.  \vspace{1mm}

Step 2 is called when an iterate $ w_j $ is obtained for which $ \lambda_{\min}(w_j) \geq 1/2 $.  Fixing such an iterate for now, let $ z := c \cdot w_j $.  $ \mathrm{AGM}_1 $ and $ \mathrm{AGM}_2 $    are started anew with input $ v_+ = e + \smfrac{3}{4} \, \smfrac{1}{1 - \lambda_{\min}(w_j)} \, ( w_j - e)  $.  Then, until Step 2 is again invoked, all iterates have objective value
 \[ 
   z_+ := c \cdot v_+ \, =  \, c \cdot e + \smfrac{3}{4} \, \smfrac{1}{1 - \lambda_{\min}(w_j)} \, ( z - c \cdot e)  \; , \]
which satisfies
\begin{equation}  \label{eqn.dh} 
   \frac{c \cdot e - z_+}{c \cdot e - z}  \,  \geq \, \frac{3}{4} \, \, \frac{1}{1 - \smfrac{1}{2} } \,  = \frac{3}{2} \; .  
   \end{equation} 
   
We observe, too, that for Step 2 to have been invoked implies 
\begin{equation}  \label{eqn.di} 
   \frac{c \cdot e - z}{ c \cdot e - z^*} \, \leq \, \frac{1}{2} \; . 
   \end{equation}
This follows since
\begin{align*}
   \frac{c \cdot e - z}{ c \cdot e - z^*} & =  1 - \frac{z - z^*}{c \cdot e - z^*} \\
                & = 1 - \lambda_{\min,z}^* \quad \textrm{(by Lemma~\ref{lem.bd})} \\
                & \leq 1 - \lambda_{\min}(w_k) \\
                & \leq 1/2 \; . 
 \end{align*}
 
Letting $ v_1, v_2, \ldots, v_{\ell'} $ be the first $ \ell' $ outer iterates computed by MainAlgo, and letting $ z_{\ell} := c \cdot v_{\ell} $, from (\ref{eqn.dh})  follows
\begin{equation}  \label{eqn.dj} 
 (c \cdot e - z_{\ell}) \geq (\smfrac{3}{2} )^{ \ell - 1} (c \cdot e - z_1 )\quad \textrm{($ \ell = 1, \ldots, \ell' $)}  \; , 
 \end{equation}
 whereas (\ref{eqn.di})  implies
 \begin{equation}  \label{eqn.dk} 
            \smfrac{1}{2} ( c \cdot e - z^*) \geq c \cdot e - z_{\ell' - 1} \; . 
 \end{equation}
 Together, (\ref{eqn.dj})  and (\ref{eqn.dk})  give  $ \smfrac{1}{2} ( c \cdot e - z^*) \geq ( \smfrac{3}{2} )^{\ell' - 2} ( c \cdot e - z_1) $, and hence,
 \begin{equation}  \label{eqn.dl} 
       \ell' \, < \, 1  +  \log_{3/2} \frac{c \cdot e - z^*}{ c \cdot e - z_1} 
       \end{equation}
(where we have used $ \log_{3/2}(1/2) < -1 $). Trivially, the number of outer iterates is thus finite.  Henceforth we use $ \ell' $ to denote the final index of the outer iterations, that is,  $ v_{\ell'} $ is the last outer iterate computed by MainAlgo.  \vspace{1mm}

Next, we claim for the decreasing sequence $ z_1, z_2, \ldots z_{\ell'} $,
\begin{equation}  \label{eqn.dm} 
    \ell \leq \ell' - 2 \quad \Rightarrow \quad      \frac{ c \cdot e - z_{\ell}}{ c \cdot e - z^*} \leq 1/3 \; . 
\end{equation}
To understand why this is true, assume $ \ell < \ell' $ and $ z_{\ell} $ does not satisfy the inequality on the right of (\ref{eqn.dm}).  Then (\ref{eqn.dh})  implies
\[ c \cdot e - z_{\ell + 1} \geq \smfrac{3}{2} ( c \cdot e - z_{\ell}) > \smfrac{1}{2} ( c \cdot e - z^*)  \; , \]
and thus, by Lemma~\ref{lem.bd},
\[  \lambda_{\min,z_{\ell + 1}}^* \, = \, 1 - \smfrac{c \cdot e - z_{\ell + 1}}{c \cdot e - z^* } \, < \, 1/2 \; . \]
Consequently, when $ \mathrm{AGM}_1 $ and $ \mathrm{AGM}_2 $    are started anew with input $ v_{\ell + 1} $, each iterate $ w_j $ must satisfy $ \lambda_{\min}(w_j) < 1/2 $, implying Step 2 will never again be invoked, that is, implying $ \ell + 1 = \ell' $.  The implication (\ref{eqn.dm})  is thus established. \vspace{1mm}

Assuming for index $ \ell $ that the inequality on the right of (\ref{eqn.dm})  is satisfied, we now provide an upper bound on the number of gradients computed by $ \mathrm{AGM}_1 $   in the $ \ell^{th} $ outer iteration, that is, the number of gradients computed between when $ \mathrm{AGM}_1 $ and $ \mathrm{AGM}_2 $  were started with input $ v_{\ell} $ and when an iterate $ w_j $ is encountered for which $ \lambda_{\min}(w_j) \geq 1/2 $. Here we rely on the equality $ \lambda_{\min}(v_{\ell})  = 1/4 $ (MainAlgo explicitly requires this equality for input $ v_1 $, and it holds for subsequent outer-iterates due to (\ref{eqn.ba})
 and to the definition of $ v_{\ell+1} $ in Step 2 of MainAlgo). \vspace{1mm}

Since $ \frac{c \cdot e - z_{\ell}}{c \cdot e - z^*} \leq \frac{1}{3}  $ (by assumption), Lemma~\ref{lem.bd}  implies
\[    \lambda_{\min,z_{\ell}}^*  = 1 - \smfrac{c \cdot e - z_{\ell}}{c \cdot e - z^*} \geq 2/3 \; . \]
Consequently, applying Corollary~\ref{cor.db}  with the choices $ x_0 = v_{\ell} $ and $ \epsilon' = 1/6 $, we find that in the $ \ell^{th} $ outer iteration, $ \mathrm{AGM}_1 $ determines an iterate satisfying $ \lambda_{\min}(w_j) \geq 1/2 $ within a number of gradient computations not exceeding 
\begin{equation} \label{eqn.dn} 
    48 \, \mathrm{diam}_{z_{\ell}} \sqrt{ 2 L_1 \ln n} + 2 \left| \log_2(  L_1^{\bullet}/L_1  ) \right| + 6 \; . 
    \end{equation}
(Upon determining $ w_j $, of course, the $ (\ell +1)^{th} $ outer iteration commences.) As $ \mathrm{AGM}_1 $ and $ \mathrm{AGM}_2 $ are assumed to compute gradients simultaneously, the same bound holds for $ \mathrm{AGM}_2 $.    \vspace{1mm}

In addition to having a bound on the number of gradients computed during outer iteration $ \ell $  when $ \frac{c \cdot e - z_{\ell}}{c \cdot e - z^*} \leq  \frac{1}{3} $, we now know for $ \ell = \ell' $ (the final outer iteration), $ \frac{c \cdot e - z_{\ell}}{c \cdot e - z^*} > \frac{1}{3} $, simply because MainAlgo never terminates.  Moreover, according to the implication (\ref{eqn.dm}), outer iteration $ \ell = \ell' - 1 $ is the only other outer iteration that might satisfy  $ \frac{c \cdot e - z_{\ell}}{c \cdot e - z^*} > \frac{1}{3} $. \vspace{1mm} 

Consider the case that $ \ell = \ell' $ is the only outer iteration for which $ \frac{c \cdot e - z_{\ell}}{c \cdot e - z^*} > \frac{1}{3} $.    Then, from above, the total number of gradients computed by $ \mathrm{AGM}_1 $   before reaching outer iteration $ \ell' $ does not exceed
\begin{equation} \label{eqn.do} 
       (\ell' - 1) \left( 48 \, \mathrm{Diam}_{z_1} \sqrt{ 2 L_1 \ln n} + 2 | \log_2(  L_1^{\bullet}/L_1  )|   + 6 \right) \; .  \end{equation}

On the other hand, since $ \frac{c \cdot e - z_{\ell'}}{ c \cdot e - z^*} > 1/3 $ and $ \lambda_{\min}(v_{\ell'}) = 1/4 > 1/6 $, Corollary~\ref{cor.dc}  implies that at some time during outer iteration $ \ell' $, an iterate $ x_k $ will be determined by $ \mathrm{AGM}_2  $ for which 
\begin{equation}  \label{eqn.dp} 
    \frac{c \cdot \pi(x_k) - z^*}{ c \cdot e - z^*} \, \leq \, \epsilon \; . 
\end{equation}
Furthermore, the corollary implies that during outer iteration $ \ell' $, the number of gradients computed by $ \mathrm{AGM}_2 $  before that time is reached is at most
\begin{equation}  \label{eqn.dq} 
24 \, \mathrm{Diam}_{z_1} \sqrt{2 L_1 \ln n}/ \epsilon  + 2 | \log_2(  L_1^{\bullet}/L_1  ) | + 6 \; .   
\end{equation}
Since $ \mathrm{AGM}_1 $ and $ \mathrm{AGM}_2 $ are assumed to compute gradients simultaneously, the same bound holds for $ \mathrm{AGM}_1 $. \vspace{1mm}

Clearly, then, for the case that $ \ell = \ell' $ is the only outer iteration for which $ \frac{c \cdot e - z_{\ell}}{c \cdot e - z^*} > \frac{1}{3}  $, the total number of gradients computed by $ \mathrm{AGM}_1 $ before (\ref{eqn.dp})  is satisfied does not exceed the sum of (\ref{eqn.do})  and (\ref{eqn.dq}).  Likewise for $ \mathrm{AGM}_2 $.  Using (\ref{eqn.dl})  to substitute for $ \ell' $ completes the proof of the theorem in this case. 
\vspace{1mm}

Finally, consider the remaining case, where $ \frac{c \cdot e - z_{\ell}}{c \cdot e - z^*} > \frac{1}{3} $ for $ \ell = \ell' - 1 $ as well as for $ \ell = \ell' $  (although $ \frac{c \cdot e - z_{\ell}}{c \cdot e - z^*} \leq  \frac{1}{3} $ for $ \ell < \ell' - 1 $).  
Clearly, replacing $ \ell' - 1 $ in (\ref{eqn.do})  by $ \ell' - 2 $ results in an upper bound on the number gradients computed by $ \mathrm{AGM}_1 $ before reaching outer iteration $ \ell' - 1 $.  \vspace{1mm}
   
If in outer iteration $ \ell = \ell' - 1 $ or $ \ell = \ell' $,  the number of gradients computed by $ \mathrm{AGM}_2 $ during the iteration reaches the quantity (\ref{eqn.dq}), then by the same reasoning as above, an iterate $ x_k $ satisfying (\ref{eqn.dp})  will have been obtained.  Consequently, the combined number of gradients computed by $ \mathrm{AGM}_2 $ in outer iterations $ \ell' - 1 $ and  $ \ell' $ before reaching an iterate $ x_k $ satisfying (\ref{eqn.dp})  is at most twice the value (\ref{eqn.dq}).  Adding this quantity to (\ref{eqn.do})  bounds the total number of gradients computed by $ \mathrm{AGM}_2 $ from the time MainAlgo is initiated until $ x_k $ satisfying (\ref{eqn.dp})  is reached.  The same bound holds for $ \mathrm{AGM}_1 $.    \hfill $ \Box $
 \vspace{5mm}

\appendix

\section{{\bf Proof of Proposition~\ref{prop.ca}}} \label{app.a}

For the proposition, $ {\mathcal K} $ is a hyperbolicity cone in $ {\mathcal E} $, with hyperbolic polynomial $ p $, of degree $ n $.  Fix $ e \in \int( {\mathcal K}) $, the distinguished direction for defining eigenvalues $ \lambda_j(x) $ for every $ x \in {\mathcal E} $, these being the roots of the univariate polynomial $ \lambda \mapsto p(x - \lambda e) $. \vspace{1mm}      

It suffices to establish the proposition for the case $ \mu = 1 $, because for any $ \mu > 0 $ we have $  f_{\mu}(x) = \mu f_1( \smfrac{1}{\mu} x) $ (due to the homogeneity $ \lambda_j( \smfrac{1}{\mu}x) = \lambda_j(x)/\mu $). Thus, we aim to show the function
\begin{equation}  \label{eqn.eaa}
  f(x) := \ln \sum_{j=1}^n \exp( \lambda_j(x)) 
  \end{equation} 
is convex and infinitely Fr\'{e}chet-differentiable, and satisfies
\begin{equation}   \label{eqn.ea} 
  \| \grad f(x) - \grad f(y) \|_{\infty}^* \leq \| x - y \|_{\infty} \quad \textrm{for all $ x, y \in {\mathcal E}$}  \; . 
  \end{equation}   
   \vspace{1mm}
 
Nesterov \cite[(10) and Thm 3]{nesterov2007smoothing}  established these facts for the special case $ {\mathcal K} = \Sym $, with hyperbolic polynomial $ X \mapsto  \det(X) $ and distinguished direction $ I $, in which case, of course, $  \lambda_j(X)   $ are the usual eigenvalues of matrix $ X $.  As we rely heavily on Nesterov's results, we capitalize the function in the special case:
\[  F(X) = \ln \sum_j \exp( \lambda_j(X)) \; . \]

Our goals are almost entirely accomplished by a straightforward leveraging of Nesterov's results using the following deep theorem. \vspace{1mm}

\begin{helton_vinnikov_thm}  
\hypertarget{targ_helton_vinnikov_thm}{} If $ {\mathcal L} $ is   a 3-dimensional subspace of $ {\mathcal E} $, and $ e \in {\mathcal L} \cap \mathrm{int}(K)   $, then there exists a linear transformation $ T: {\mathcal L}  \rightarrow \Sym $ satisfying 
 \[     T(e) = I \quad \textrm{and} \quad p(x) = p(e) \, \det(T(x)) \, \textrm{ for all $ x \in {\mathcal L}  $} \; . \]
\end{helton_vinnikov_thm}
\vspace{3mm}

\addtocounter{prop}{1}
 
This is a ``homogeneous version'' of the theorem  proven by Helton and Vinnikov \cite{helton2007linear} (also see \cite{vinnikov2012lmi}).  The homogenous version was first recorded by Lewis, Parrilo and Ramana in \cite{lewis2005lax}, and observed there to settle affirmatively the Lax conjecture (i.e., Peter Lax had conjectured that the statement of the above theorem is true). \vspace{1mm}

Recently, Hanselka \cite{hanselka2014definite}  gave a proof of the H-V (Helton-Vinnikov) Theorem very different than the original proof, his approach being largely algebraic rather than geometric.  Both proofs rely on deep mathematics.  \vspace{2mm}

A simple observation is that convexity is a 1-dimensional property -- a function with domain $ {\mathcal E} $ is convex if and only if it is convex when restricted to lines in $ {\mathcal E} $.  Thus, to prove convexity of $ f $, it suffices to consider general pairs $ x, y \in {\mathcal E} $ and show $ f $ is convex when restricted to the line containing $ x $ and $ y $, which certainly is the case if $ f $ is convex on the subspace $ {\mathcal L} $  spanned by $ x, y $ and $ e $.  However, letting $ T $ be as in the H-V Theorem, it is easy to see $ \{ \lambda_j(x) \}  = \{ \lambda_j(T(x)) \}  $ for every $ x \in {\mathcal L} $, and thus, $ f(x) = F(T(x)) $. Convexity of $ f $ on $ {\mathcal L} $ is then immediate from convexity of $ F $.\footnote{Invoking the Helton-Vinnikov Theorem to establish convexity of $ f $ is overkill.  For a more insightful understanding of the convexity of $ f $ (and the convexity of a vast range of other ``eigenvalue'' functions), see \cite{bauschke2001hyperbolic}.}  \vspace{1mm}

Likewise, for any convex function $ f \in  {\mathcal C}^1 $, Lipschitz continuity of the gradient map is a 1-dimensional property, in that, for example, for any norm $ \| \, \, \|_{\vartriangle} $ and its dual $ \| g \|_{\vartriangle}^* := \max \{ u: \lin g, u \rin \leq 1 \} $,  
\begin{gather}  \| \grad f(x) - \grad f(y) \|_{\vartriangle}^* \leq L \| x - y \|_{\vartriangle} \quad \textrm{for all $ x, y \in {\mathcal E} $}  \nonumber \\
\Leftrightarrow  \nonumber \\
 f(y) \leq f(x) + \smfrac{d}{dt} f(x+t(y-x))|_{t=0} + \smfrac{L}{2} \| y - x \|_{\vartriangle}^2 \quad \textrm{for all $ x, y \in {\mathcal E}$} \; .  \label{eqn.ec} 
 \end{gather}
(The implication ``$ \Rightarrow $'' is standard, but the implication ``$ \Leftarrow $'' holds as well, as is shown in \cite{nesterov2004introductory}  for Euclidean norms, and is remarked in \cite[Appendix]{nesterov2007smoothing}  to hold generally\footnote{Indeed, in the proof of \cite[Thm. 2.1.5]{nesterov2004introductory}, for the function $ \phi(y) := f(y) - \lin \grad f(x), y - x \rin $, simply replace the point $ y - \smfrac{1}{L} \grad \phi(y) $ with $ y - \smfrac{\| \grad \phi(y)\|_{\vartriangle}^*}{L} u $ where $ \| u \|_{\vartriangle} = 1 $ and $ \lin \grad \phi(y), u \rin = \| \grad \phi(y) \|_{\vartriangle}^* $.}.) \vspace{1mm}

Thus, to establish the desired Lipschitz continuity (\ref{eqn.ea})  for the gradient map of the function $ f $ defined by (\ref{eqn.eaa}) , it suffices to show $ f \in {\mathcal C}^1 $ and to establish the inequality (\ref{eqn.ec})  for a general pair $ x, y \in {\mathcal E}$, with $ L = 1 $ and $ \| \, \, \|_{\vartriangle} = \| \, \, \|_{\infty} $. \vspace{1mm}

However, assuming for the moment that $ f $ is continuously differentiable, letting $ {\mathcal L} $ be the subspace spanned by $ x, y $ and $ e $, and letting $ T $ be as in the H-V Theorem, for $ X := T(x) $ and $ Y := T(y) $ we have 
\begin{align*}
   f(y) & = F(Y) \\
          & \leq F(X) + \smfrac{d}{dt} F(X+t(Y-X))|_{t=0} + \smfrac{1}{2} \| Y - X \|_{\infty}^2 \\  & = f(x) + \smfrac{d}{dt} f(x+t(y-x))|_{t=0} + \smfrac{1}{2} \| y - x \|_{\infty}^2 \; , 
          \end{align*}
where the inequality is due to Nesterov \cite[(10) and Thm 3]{nesterov2007smoothing}  already having established the desired Lipschitz continuity for the gradient map of $ F $. Hence, to complete the proof of the desired Lipschitz continuity for gradient map of $ f $, it only remains to show $ f \in {\mathcal C}^1 $.   \vspace{1mm}

We show that $ f $ is, in fact, analytic.\footnote{In the original version we showed only that $ f \in {\mathcal C}^{\infty} $. We thank a referee for observing that the proof could be both streamlined and made to lead to the stronger result of analyticity.}  For this it clearly suffices to show $ x \mapsto \sum_j \exp(\lambda_j(x)) $ is analytic, that is, it suffices to establish analyticity of 
\begin{equation}  \label{eqn.azaz} 
x \mapsto \sum_{k=0}^{\infty} \smfrac{1}{k!} S_k(x) \quad \textrm{where } S_k(x) := \sum_j \lambda_j(x)^k \; . \end{equation}  

For establishing analyticity of (\ref{eqn.azaz}), we first note it is well known that the function $ S_k $ is a homogeneous polynomial of degree $ k $  (as follows, for example, from Fact 2.10 and Theorem 3.1 in \cite{bauschke2001hyperbolic}).  Hence, the functions
\[ 
  T_{\ell}(x) :=  \sum_{k=0}^{\ell} \smfrac{1}{k!} S_k(x) \quad  \textrm{($ \ell = 0, 1, \ldots $)} \; 
\] 
are polynomials, and thus extend to arguments $ x \in \mathbb{C}^n $.  To establish analyticity of the function (\ref{eqn.azaz}), it suffices to show the holomorphic functions $ T_{\ell} $ converge uniformly on compact subsets of $ \mathbb{C}^n $, because it is well known that the limit function will then be holomorphic (c.f., \cite[Prop 2]{malgrange1984lectures}), and hence analytic when restricted to $ \mathbb{R}^n $.

To this end, observe
\[   M :=   \sup \{  \sum_j | \lambda_j(x)| : x \in \mathbb{C}^n \textrm{ and } \sum_i | x_i | \leq 1 \} \, < \, \infty,      \]
finiteness due to compactness of the set of arguments $ x $, and the fact that roots of a 
univariate polynomial vary continuously in the coefficients if the leading coefficient remains constant -- and thus the roots of $ \lambda \mapsto p(x - \lambda e) $ vary continuously in $ x $. 
Clearly, for every $ k \geq 1 $,
\[  
    \sup \{  |S_k(x)| : x \in \mathbb{C}^n \textrm{ and } \sum_i |x_i| \leq 1 \}    \, \leq \, M^k  \; . 
\] 
Consequently, as $ S_k $ is homogeneous of degree $ k $, 
\[   | S_k(x) | \leq \left(  M \sum_i |x_i | \right)^k \quad \textrm{for all $ x \in \mathbb{C}^n $} \; . \]
It is readily argued from this and the definition of $ T_{ \ell} $ that the sequence $ \{ T_{ \ell } \} $ converges uniformly on compact subsets of $ \mathbb{C}^n $, completing the proof.  \hfill $ \Box $
\vspace{5mm}

\section{{\bf  Proof of Propositon~\ref{prop.cc}}}   \label{app.b}
For ease of reference, we restate the proposition.

\begin{prop} For $ x \in {\mathcal E} $, assume $ \{ \lambda_j(x) \} $ is the set of \underline{distinct} eigenvalues of $ x $, and let $ m_j $ denote the multiplicity of $ \lambda_j(x) $ (that is, the multiplicity of $ \lambda_j(x) $ as a root of $ \lambda \mapsto p(x - \lambda e) $).  Then
\begin{equation}  \label{eqn.faa}
    \grad f_{\mu}(x) = \smfrac{1}{\sum_j m_j \exp(\lambda_j(x)/\mu)} \sum_j \smfrac{m_j \exp(\lambda_j(x)/\mu)}{p^{(m_j)}( \, x - \lambda_j(x) e \, )} \, \grad p^{(m_j - 1)}( \, x - \lambda_j(x) e \, ) \; , 
    \end{equation} 
a convex combination of the vectors $ \smfrac{1}{p^{(m_j)}( \, x - \lambda_j(x) e \, )} \, \grad p^{(m_j - 1)}( \, x - \lambda_j(x) e \, ) $ (vectors which are independent of $ \mu $).
\end{prop}
\noindent {\bf Proof:}  We first establish that the statement is correct for the hyperbolic polynomial $ P(X) = \det(X) $ ($ X \in \Sym $) and distinguished direction $ I $, and then rely on the Helton-Vinnikov Theorem (Appendix~\ref{app.a}) to establish the theorem generally. \vspace{1mm}

For $ P(X) = \det(X) $, it is readily verified by induction that the value $ P^{(i)}(X) $ is $ i! $ ($ i $ factorial) times the sum of the $ (n-i)^{th} $ order principal minors of $ X $.  Thus, assuming $ D $ is a diagonal matrix having zero for its first $ m $ diagonal entries, and having nonzero eigenvalues $ \lambda_{m+1}, \ldots, \lambda_n $, we have $ P^{(m)}(D) = m! \, \prod_{i=m+1}^n \lambda_i $.  Moreover, for gradients computed with respect to the trace inner product, it is straightforward (even if notationally cumbersome) to show $ \grad P^{(m-1)}(D) $ is the diagonal matrix whose first $ m $ diagonal entries are $ (m-1)! \prod_{i=m+1}^n \lambda_i $ and whose other entries are all zero. Consequently, $   
\smfrac{m}{P^{(m)}(D)} \grad P^{(m-1)}(D) $ is the diagonal matrix whose first $ m $ diagonal entries are 1, and whose other entries are 0. \vspace{1mm}

For any $ n \times n $ orthogonal matrix $ Q $, it follows that 
\[ \smfrac{m}{P^{(m)}(Q D Q^T)} \grad P^{(m-1)}(Q D Q^T) = Q_{\cdot ,1:m} (Q_{\cdot ,1:m})^T \; , \]
where $ Q_{\cdot ,1:m} $ is the matrix formed from the first $ m $ columns of $ Q $.  Thus,  \\
 $ \smfrac{m}{P^{(m)}(Q D Q^T)} \grad P^{(m-1)}(Q D Q^T) $ is the matrix projecting $ \mathbb{R}^n $ orthogonally onto the null space of $ QD Q^T $. \vspace{1mm}

 Consequently, if for $ X \in \Sym $, $ \lambda_j(X) $ is an eigenvalue of multiplicity $ m_j $, then the matrix $ \smfrac{m_j}{P^{(m_j)}(X - \lambda_j(X)I )} \grad P^{(m_j-1)}(X - \lambda_j(X) I )  $ projects $ \mathbb{R}^n $ orthogonally onto the eigenspace corresponding to $ \lambda_j(X) $.  \vspace{1mm}

Since for the function $ F_{\mu}(X) = \mu \ln \sum_j \exp(\lambda_j(X)/\mu) $, it already is known that 
\[   \grad F_{\mu}(X) = \smfrac{1}{\sum_j \exp(\lambda_j(X)/\mu)} Q \left[ \begin{smallmatrix} \exp(\lambda_1(X)/\mu) & & \\ & \ddots & \\ && \exp(\lambda_n(X)/\mu) \end{smallmatrix} \right]   Q^T \]  
(where $ Q \left[ \begin{smallmatrix} \lambda_1(X) & & \\ & \ddots & \\ & & \lambda_n(X) \end{smallmatrix} \right]  Q^T $ is an eigendecomposition of $ X $), the proposition immediately follows for the special case of the hyperbolic polynomial $ P(X) = \det(X) $ and distinguished direction $ I $. \vspace{1mm}
  
More generally, let $ {\mathcal K}    $ be a hyperbolicity cone in a Euclidean space $ {\mathcal E} $, and assume $ p $ is a hyperbolic polynomial for $ {\mathcal K}  $.  Let $ n $ be the degree of $ p $, and assume $ e \in \int({\mathcal K}) $, the distinguished direction for defining eigenvalues. \vspace{1mm}

We now briefly argue that in proving the proposition, we may assume $ {\mathcal K} $ is a regular cone, that is, contains no subspace other than $ \{ \vec{0} \} $ and has nonempty interior (the latter condition being satisfied simply because $ e \in \int({\mathcal K}) $, by assumption). \vspace{1mm}

Let $ {\mathcal L} $ denote the lineality space of $ {\mathcal K} $, the largest subspace contained in $ {\mathcal K} $.  Then for every $ e' \in e + {\mathcal L} $ (Minkowski sum), and for every $ x \in {\mathcal E} $, the eigenvalues for $ x $ as defined by $ e' $ are the same as those defined by $ e $.  Moreover, for every $ x' \in x + {\mathcal L} $, the eigenvalues of $ x' $ are the same as those of $ x $. (For simple justification of these claims, see \cite[\S 8]{renegar2014polynomial}.) \vspace{1mm}

Since the eigenvalues defined by $ e' \in e + {\mathcal L}   $ are the same as those defined by $ e $, also the same are the values of the function $ f_{\mu} $ and polynomials $ p^{(i)} $.  Consequently, we may assume $ e \in {\mathcal L}^{\perp} $ (where orthogonality is with respect any fixed inner product). \vspace{1mm}

Since the eigenvalues of $ x' \in x + {\mathcal L}  $ are the same as those of $ x $, and since this holds for all $ x \in {\mathcal E} $, we have $ \grad f_{\mu}(x') = \grad f_{\mu}(x) $ and $ \grad p^{(i)}(x') = \grad p^{(i)}(x) $ for all $ x' \in x + {\mathcal L} $, and we have that the gradients lie in $ {\mathcal L}^{\perp} $. \vspace{1mm}

In light of these observations, in proving the proposition we may restrict attention to $ x $ lying in the cone $ {\mathcal K} \cap {\mathcal L}^{\perp} $, a regular cone in the Euclidean space $ {\mathcal E} \cap {\mathcal L}^{\perp} $, and which has hyperbolic polynomial $ p|_{{\mathcal E} \cap {\mathcal L}^{\perp}} $.  Thus, as claimed, we may assume $ {\mathcal K} $ is a regular cone. \vspace{1mm}

Although gradients depend on the particular inner product, the desired identity (\ref{eqn.faa})  does not.  This allows us to proceed with an inner product of our choosing.  We choose the inner product $ \lin u, v \rin_e := \lin u, H(e) v \rin $, where $ H(e) $ is the Hessian at $ e $ for the self-concordant logarithmic barrier function $ x \mapsto - \ln p(x) $.  (Regularity of $ {\mathcal K} $ implies invertibility of $ H(e) $, and thus guarantees $ \langle \; , \; \rangle_e $ is indeed an inner product -- elementary proofs of these well-known facts can be found in \cite[\S 8]{renegar2014polynomial}.) \vspace{1mm} 

Fix $ x \in {\mathcal E} $, and let 
\[ g_1 := \grad f_{\mu}(x) \; , \quad  g_2 :=  \smfrac{1}{\sum_j m_j \exp(\lambda_j(x)/\mu)} \sum_j \smfrac{m_j \exp(\lambda_j(x)/\mu)}{p^{(m_j)}( \, x - \lambda_j(x) e \, )} \, \grad p^{(m_j - 1)}( \, x - \lambda_j(x) e \, ), \]
where $ \{ \lambda_j(x) \} $ are the distinct eigenvalues of $ x $, and where $ m_j $ is the multiplicity of $ \lambda_j(x) $.  Our goal is to show $ g_1 = g_2 $. \vspace{1mm}

Consider a subspace $ {\mathcal L} $ whose dimension is at most 3, and which satisfies $ e,x \in {\mathcal L} \subseteq  {\mathcal E} $.   Let $ T: {\mathcal L} \rightarrow \mathbb{R} $ be a linear transformation as in the Helton-Vinnikov Theorem:  
\begin{equation} \label{eqn.fa} 
  p(e) = I \quad \textrm{and} \quad   p(y) = p(e) \det( T(y)) \textrm{ for all $ y \in {\mathcal L} $} \; . \end{equation}
  
  From (\ref{eqn.fa})  and the definition of $ \langle \; , \; \rangle_e $ follows
\[  
   \lin u, v \rin_e  = \textrm{trace}(T(u) T(v)) \, \textrm{ for all $ u, v \in {\mathcal L} $}, \] 
that is, $ T $ is an isometry between $ {\mathcal L} $ (endowed with the inner product $ \langle \; , \; \rangle_e $) and the subspace $ T( {\mathcal L}) $ (endowed with the trace inner product).  Consequently, if $ {\mathcal C}^1 $   functions $ h: {\mathcal E} \rightarrow \mathbb{R} $ and $ H: \Sym \rightarrow \mathbb{R} $ satisfy $ h(y) = H(T(y)) $ for all $ y \in {\mathcal L} $, then for all $ y \in {\mathcal L} $, 
\begin{equation} \label{eqn.fb} 
  T( \textrm{proj}_{{\mathcal L}} \grad h(y)) = \textrm{Proj}_{T( {\mathcal L})} \grad H(T(y)) \; , 
  \end{equation} 
where $ \textrm{proj}_{{\mathcal L}} $ (resp., $ \textrm{Proj}_{T({\mathcal L})} $)  is the operator projecting $ {\mathcal E} $ orthogonally onto $ {\mathcal L} $ (resp., projecting $ \Sym $ orthogonally onto $ T({\mathcal L}) $).   \vspace{1mm}

On the other hand, clear from (\ref{eqn.fa})  is that for every $ i $,
 \begin{equation}  \label{eqn.fc}  
   p^{(i)}(y) = p(e) P^{(i)}(T(y)) \, \textrm{ for all $ y \in {\mathcal L} $} \; . 
   \end{equation}
 Also evident is  
 \begin{equation}  \label{eqn.fd} 
   f_{\mu}(y) = F_{\mu}(T(y)) \, \textrm{ for all $ y \in {\mathcal L} $} \; . 
   \end{equation}

Letting $ X := T(x) $ and $ G_1 := \grad F_{\mu}(X) $, from (\ref{eqn.fb})  and (\ref{eqn.fd})  follows
\[  
 T( \textrm{proj}_{{\mathcal L}} g_1 ) = \textrm{Proj}_{T({\mathcal L})} G_1 \; . \]
Similarly, from (\ref{eqn.fb})  and (\ref{eqn.fc})  follows 
\[  
 T( \textrm{proj}_{{\mathcal L}} g_2 ) = \textrm{Proj}_{T({\mathcal L})} G_2 \; , \]
where 
\[ G_2 := \smfrac{1}{\sum_j m_j \exp(\lambda_j(X)/\mu)} \sum_j \smfrac{m_j \exp(\lambda_j(X))/\mu)}{P^{(m_j)}( \, X - \lambda_j(X) I \, )} \, \grad P^{(m_j - 1)}( \, X - \lambda_j(X) I \, ) \; .
\] 
However, we have already seen $ G_1 = G_2 $.  Thus,  for any $ 3 $-dimensional subspace $ {\mathcal L} $ containing $ e $ and $ x $, 
\begin{equation}  \label{eqn.fe} 
  \textrm{proj}_{{\mathcal L}} g_1 = \textrm{proj}_{{\mathcal L}} g_2 \; . 
  \end{equation}
  
Letting $ {\mathcal L} $ be the subspace spanned by $ e $, $ x $ and $ g_1 $, it follows from (\ref{eqn.fe})  that $ g_1 $ is the orthogonal projection of $ g_2 $ onto the line $ \{ tg_1: t \in \mathbb{R} \} $.  Similarly, letting $ {\mathcal L} $ be the subspace spanned by $ e $, $ x $ and $ g_2 $, we find $ g_2 $ is the orthogonal projection of $ g_1 $ onto the line $ \{ t g_2: t \in \mathbb{R} \} $. We conclude $ g_1 = g_2 $. \hfill $ \Box $
\vspace{5mm}

\section{{\bf  Proof of Theorem~\ref{thm.da}}}  \label{app.c}

Theorem~\ref{thm.da}, due to Nesterov, is not quite to be found in \cite{nesterov2013gradient}  although he develops all of the key ingredients.  He assumes, for example, that the estimated value for the Lipschitz constant is a lower bound, but does not record the (easier) case that the estimated value is an upper bound.  \vspace{1mm}

Additionally, the development in \cite{nesterov2013gradient}  is done in the context of general prox functions, whereas we need only the prox function $ x \mapsto \smfrac{1}{2} \| \, \, \|^2 $, as our problems have been recast to have only linear equations as constraints and hence nothing more than orthogonal projection of gradients onto a subspace is needed to ensure iterates remain feasible. \vspace{1mm}

Because Theorem~\ref{thm.da}  is not quite to be found in \cite{nesterov2013gradient}, and because the development there can be streamlined in our setting, in this appendix we provide a proof of the theorem, relying solely on Nesterov's ideas.  \vspace{1mm}

Let $ {\mathcal E} $ be a Euclidean space, with inner product $ \langle \; , \; \rangle $ and associated norm $ \| \, \, \| $.  Let $ f : {\mathcal E} \rightarrow \mathbb{R} $ be a \underline{convex} function with Lipschitz-continuous gradient.  Let $ L_f $ be the Lipschitz constant (possibly unknown to the user).  The goal is simply to minimize $ f $  with no constraints (linear equations are easily handled in the theory by taking projected gradient steps rather than gradient steps). \vspace{1mm}

Following is the accelerated method \cite[(4.9)]{nesterov2013gradient}, streamlined for this simple setting.

\noindent 
\hrulefill
 
\noindent {\bf Accelerated Gradient Method} \vspace{1mm}

\noindent  
(0) Input: $ x_0 \in {\mathcal E} $ and $ L' > 0 $, an estimate of $ L_f $.  \\
$ \textrm{~} $ \quad Initialize: $ v = x_0 $, \, $ L = L' $,  \, $ A = 0 $ \,
and \,
 $ k = 0 \; .   $    \\
(1) Compute $ a = \frac{1 + \sqrt{1 + 2 AL}}{L} $, \,
 $ y = \frac{A}{A + a} x_k + \frac{a}{A + a} v $ \,
  and \,
   $ x = y - \smfrac{1}{L}  \grad f(y) \; . $    \\  
(2) If $ \lin \grad f(x), y - x \rin < \frac{1}{L} \,  \| \grad f(x) \|^2 $   then let $ L \leftarrow 2L $ and return to Step 1. \\
(3) Define $ x_{k+1} = x  $ and let \, $ v \leftarrow v - a \, \grad f(x)  $,  \, $ L \leftarrow L/2 $, \, $ A \leftarrow A + a \; , $     \\
$ \textrm{~} $ \qquad \qquad  \qquad  \qquad  \qquad  \qquad  \qquad  \qquad  \qquad  \qquad  \qquad  $ k \leftarrow k+1 \; , $ and go to Step 1.   
  \vspace{-1.5mm}

\noindent 
\hrulefill
\vspace{2mm}

\noindent {\bf Remark:} Whereas in Step 2 we have $ L \leftarrow 2L $, Nesterov has $ L \leftarrow \gamma_u L $, where $ \gamma_u > 1 $ is user-chosen.  Similarly, in Step 3 he has $ L \leftarrow L/ \gamma_d $.  He assumes $ \gamma_d \geq \gamma_u $. \vspace{1mm}

\begin{thm}  {\bf  (Nesterov \cite{nesterov2013gradient})}   \label{thm.ga}
Let  $ k' := \max \{ 0, \lceil  \log_2 (L'/L_f) \rceil \}  $. Then 
\begin{equation}  \label{eqn.ga}
   k > k' \quad \Rightarrow \quad   f(x_k) - f^* \, \leq \, 2 L_f  \, \left(  \frac{ \| x_0 - x^* \|}{ k-k'} \right)^2  \; . 
   \end{equation} 
Moreover, for $ k > k' $, the total number of gradients computed before $ x_k $ is determined does not exceed 
\begin{equation}  \label{eqn.gb}
   4 (k - k') + 2 \,  | \log_2(  L'/L_f |  ) | + 6 \; . 
   \end{equation}  
\end{thm}
\vspace{2mm}

In establishing the theorem, we rely on Nesterov's notation, as given in the following rendition of the algorithm. 
 
\noindent 
\hrulefill
 
\noindent {\em (Notationally-embellished version of the accelerated method)}  \vspace{1mm}

\noindent  
(0) Input: $ x_0 \in {\mathcal E} $ and $ L' > 0 $, a guessed estimate of $ L_f $.  \\
$ \textrm{~} $ \quad Initialize: $ v_0 = x_0 $, \, $ L_0 = L' $, \,  $ L = L_0 $, \, $ A_0 = 0 $, \, 
 $ k = 0    $.   \\
(1) Compute $ a = \frac{1 + \sqrt{1 + 2 A_k L}}{L} $, \,
 $ y = \frac{A_k}{A_k + a} x_k + \frac{a}{A_k + a} v_k $ \,   and \,
   $ x = y - \smfrac{1}{L}  \grad f(y) \; . $    \\  
(2) If $ \lin \grad f(x), y - x \rin < \frac{1}{L} \,  \| \grad f(x) \|^2 $   then let $ L \leftarrow 2L $ and return to Step 1. \\
(3) Define $ y_k = y $, \, $ M_k = L $, \\
$ \textrm{~} $ \qquad  $ x_{k+1} = x  $,  \, $ v_{k+1} = v_k - a \, \grad f(x)  $, \,  $ L_{k+1} = L/2 $, \, $ a_{k+1} = a $, \, $ A_{k+1} = A_k + a $,  \\
 $ \textrm{~} $ \quad \,  and let $ L \leftarrow L/2 $, $ k \leftarrow k+1 $, then go to Step 1.   
  \vspace{-1.5mm}

\noindent 
\hrulefill
\vspace{2mm}

Inductively define the convex functions $ \psi_k $ by $ \psi_0(x) := \smfrac{1}{2} \| x - x_0 \|^2 $ and 
\[  \psi_{k+1}(x) := \psi_k(x) + a_{k+1} \big( f(x_{k+1}) +  \lin \grad f(x_{k+1}), x - x_{k+1} \rin \big) \; . \]
It is easy to verify inductively that $ v_k $ is the minimizer of $ \psi_k \; . $  \vspace{1mm}

Note that $ x \mapsto \frac{1}{A_{k+1}} \psi_{k+1}(x) $ can be thought of as an approximation to $ f $, in that
\[  \smfrac{1}{A_{k+1}}  \psi_{k+1}(x) = \smfrac{1}{2A_{k+1}} \| x - x_0 \|^2 + \smfrac{1}{\sum_{j=1}^{k+1} a_j} \sum_{j=1}^{k+1} a_j \big( f(x_j) + \lin \grad f(x_j), x - x_j \rin \big) \; , \]
a quadratic function plus a convex combination of linear approximations to $ f $ (the weight $ 1/A_{k+1} $ on the quadratic term, as will be seen, goes to zero as $ k \rightarrow \infty $).  \vspace{1mm}

 The proof of the implication (\ref{eqn.ga})  is accomplished by inductively establishing the following inequalities and implication:
\begin{gather} 
      \psi_k(x) \leq A_k f(x) + \smfrac{1}{2} \| x - x_0 \|^2 \quad \textrm{for all $ x \in {\mathcal E} $} \; ,  \label{eqn.gc} \\
           A_k f(x_k) \leq \psi_k^* \; , \label{eqn.gd} \\
  k > k' \quad \Rightarrow \quad A_k \geq \smfrac{1}{4} (k - k')^2/L_f \; ,   \label{eqn.ge}
\end{gather}  
where $ \psi_k^* $ is the optimal value of $ \psi_k $, and where $ k' $ as in the theorem, that is, $ k' =  \max \{ 0, \lceil  \log_2 (L'/L_f) \rceil \}  $. 
Clearly, (\ref{eqn.gc})  and (\ref{eqn.gd})   imply $ A_k f(x_k) \leq A_k f^* + \smfrac{1}{2} \| x^* - x_0 \|^2 $, which combined with (\ref{eqn.ge})  gives (\ref{eqn.ga}). \vspace{1mm}

\begin{prop}  Both (\ref{eqn.gc})  and (\ref{eqn.gd})  are valid. 
\end{prop}
\noindent {\bf Proof:}  To establish (\ref{eqn.gc}), first note the inequality trivially holds for $ k = 0 $.  Assuming the inequality holds for $ k $, we have for all $ x $,
\begin{align*}
   \psi_{k+1}(x) & = \psi_k(x) + a_{k+1} ( f(x_{k+1}) + \lin \grad f(x_{k+1}, x - x_{k+1} \rin ) \\
                 & \leq \psi_k(x) + a_{k+1} f(x) \quad \mathrm{(convexity)}   \\
                 & \leq (A_k + a_{k+1}) f(x) + \smfrac{1}{2} \| x - x_0 \|^2 \\
                  & = A_{k+1} f(x) +  \smfrac{1}{2} \| x - x_0 \|^2 \; . 
\end{align*}

Now we establish the inequality (\ref{eqn.gd}), which trivially holds as equality when $ k = 0 $. Assume the inequality holds for $ k $.                                           
To ease notation, remove the subscript ``$ k $'' and use ``$ + $'' in place of ``$ k+1 $.'' (Thus, for example, $ \psi_k $ becomes $ \psi $, and $ \psi_{k+1} $ becomes $ \psi_+ $). \vspace{1mm}

  Since $ v_+ $ minimizes $ \psi_+ $, 
\begin{equation} \label{eqn.gf}
   \psi_+^* = \psi(v_+) + a_+ \big( f(x_+) + \lin \grad f(x_+), v_+ - x_+ \rin \big) \; . 
   \end{equation} 
However, since $ \psi $ is the sum of a linear function and the quadratic $ x \mapsto \smfrac{1}{2} \| x - x_0 \|^2 $, and since $ v $ minimizes $ \psi $,
\begin{align*}
  \psi(v_+) & = \psi^* + \smfrac{1}{2} \| v_+ - v \|^2 \\
                  & \geq A f(x) +  \smfrac{1}{2} \| v_+ - v \|^2 \quad \textrm{(inductive hypothesis)} \\
                  &  \geq  A \big( f(x_+) + \lin \grad f(x_+), x - x_+ \rin \big) + \smfrac{1}{2} \| v_+ - v \|^2 \quad \textrm{(convexity)} \; , 
 \end{align*}                  
and hence (\ref{eqn.gf})  implies (using $ A + a_+ = A_+ $) 
\begin{align*}
   \psi_+^* & \geq A_+ f(x_+) +   \smfrac{1}{2} \| v_+ - v \|^2  \\
  &  \qquad + \lin \grad f(x_+), A ( x - x_+) + a_+ ( v_+ - x_+) \rin \\
   & = A_+ f(x_+)  \\
   &  \qquad  +  A_+ \lin \grad f(x_+),  y - x_+) \rin - \smfrac{1}{2} a_+^2 \| \grad f(x_+) \|^2 \\
& \geq  A_+ f(x_+)  \; , 
  \end{align*}
where the equality is obtained by substituting $ v_+ = v - a_+ \grad f(x_+) $ and then
\[ 
   A( x - x_+) + a_+( v - x_+) \,  = \, A x + a_+ v - A_+ x_+ 
                                    \, = \, A_+ (y - x_+) \; , \] 
  and where the final inequality comes from the condition in Step 2 allowing passage to Step 3:
\[ 
  \lin \grad f(x_+), y - x_+ \rin \,  \geq \, \smfrac{1}{M_+} \| \grad f(x_+) \|^2 \,  = \, \smfrac{ a_+^2 }{2A_+} \| \grad f(x_+) \|^2 \; . 
                                           \]                                    
The proof is complete. \hfill $ \Box $
 \vspace{5mm}
 
 To establish the implication (\ref{eqn.ga}), it remains to establish (\ref{eqn.ge}). This is done using \cite[Lemma 5]{nesterov2013gradient}, which shows that if $ L \geq L_f$, then the algorithm moves immediately from Step 2 to Step 3.  Here is a proof of the lemma in our simple setting. 

\begin{lemma} \label{lem.gc} Let $ L > 0 $, $ y \in {\mathcal E} $ and $ x = y - \smfrac{1}{L} \grad f(y) $. Then
\[   L \geq L_f \quad \Rightarrow \quad  \lin \grad f(x), y - x \rin \geq \frac{1}{L} \,  \| \grad f(x) \|^2 \; . \] 
\end{lemma}
\noindent {\bf Proof:}  The proof makes use of the fact that for all $ x, y \in {\mathcal E} $ and $ L \geq L_f $, 
\[  
        \lin \grad f(y) - \grad f(x), y - x \rin \geq \smfrac{1}{L}  \| \grad f(x) - \grad f(y) \|^2 
\] 
(see Theorem 2.1.5 in Nesterov's introductory book \cite{nesterov2004introductory}).  Substituting $ y - x = \smfrac{1}{L} \grad f(y) $, multiplying both sides by $ L $  and then rearranging, we have
\[ \lin \grad f(x), \grad f(y) \rin \geq \| \grad f(x) \|^2 \; . \]
Substituting $ \grad f(y) = L (y - x) $ completes the proof. \hfill $ \Box $
 \vspace{2mm}
 
 Recall $ k' := \max \{ 0, \lceil \log_2( L'/ L_f) \rceil  \} = \max \{ 0, \lceil \log_2( L_0/ L_f) \rceil  \} $. \vspace{1mm}

\begin{lemma} \label{lem.gd}
  The first index $ k $ satisfying $ L_k \leq L_f $ is $ k = k' $.  Moreover, for all $ k \geq  k' $, both $ L_k \leq L_f $ and $ M_k < 2 L_f $. 
\end{lemma}
\noindent {\bf Proof:}   
Lemma \ref{lem.gc}  implies that for $ L_k \geq L_f $, the algorithm goes immediately from Step 2 to Step 3, that is, $ L $ is not replaced by $ 2L $ but instead $ L_{k+1} $ is defined as $ L_k/2 $.  It follows that $ k' $ is the first index $ k $ for which $ L_k \leq  L_f $.  \vspace{1mm}

To establish the claimed bounds on $ L_k $ and $ M_k $ when $ k \geq k' $, we proceed by induction, assuming $ L_k \leq L_f $ and showing it follows that $ M_k <  2 L_f $, and hence, $ L_{k+1} = M_k/2 < L_f $. \vspace{1mm}

Assume for index $ k $ that $ L_k \leq L_f $.  Let $ m_k $ be the non-negative integer satisfying $ M_k = 2^{m_k} L_k $. Thus, $ m_k $ is the number of times of times that beginning with $ L = L_k $, the value $ L $ was multiplied by 2 before the inequality in Step 2 no longer was satisfied.

If $ m_k = 0 $, then trivially, $ M_k = L_k \leq L_f < 2 L_f $, the desired inequality. On the other hand, if $ m_k > 1 $, then $ 2^{m_k - 1} L_k < L_f $, since otherwise Lemma \ref{lem.gc} implies the inequality in Step 2 would not have been satisfied when $ L = 2^{m_k - 1} L_k $, at which time the algorithm would have preceded to Step 3 rather than again multiplying by 2. \hfill $ \Box $
 \vspace{2mm}

Now we can establish (\ref{eqn.ge}), and thus complete the proof of the implication (\ref{eqn.ga}) in Theorem~\ref{thm.ga}.
\begin{lemma} 
The implication (\ref{eqn.ge})  holds, that is,
\[  k > k' \quad \Rightarrow \quad A_k \geq \frac{(k - k')^2}{4 L_f} \; . \]
\end{lemma}
\noindent {\bf Proof:}  In the  proof of \cite[Lemma 8]{nesterov2013gradient}, Nesterov easily shows that regardless of the initial value $ L_0 $, 
\[    1 \leq 2 M_k \left( A_{k+1}^{1/2} - A_k^{1/2} \right)^2 \; .  \]
Since $ M_k \leq 2 L_f $ whenever $ k > k' $ (Lemma~\ref{lem.gd}), and since $ A_{k+1} > A_k $, we thus have
\[  k > k' \quad \Rightarrow \quad     A_{k+1}^{1/2} \geq A_k^{1/2} + \frac{1}{2 \sqrt{L_f}} \; , \]
from which the lemma follows by induction.     \hfill $ \Box $
 \vspace{1mm}

To complete the proof of Theorem~\ref{thm.ga}, it remains to prove the bound (\ref{eqn.gb})  on the number of gradient evaluations.  Following Nesterov's choice of notation, let $ N_k $ denote the total number of gradients computed before $ x_{k+1} $ becomes defined.  The following lemma establishes the desired bound. 

\begin{lemma} 
\[ 
        k \geq k' \quad \Rightarrow \quad N_k < 4(k - k')  + 2 \, |\log_2 (  L_0/L_f ) |  + 6 \] 
\end{lemma}
\noindent {\bf Proof:}    
For the case $ L_0 \leq L_f $, the proof of \cite[Lemma 6]{nesterov2013gradient}  establishes a slightly stronger result (choose $ \gamma_u = \gamma_d = 2 $ as the parameter values there).  Thus, it remains to establish the lemma when $ L_0 > L_f $. \vspace{1mm}

In the proof of \cite[Lemma 6]{nesterov2013gradient}, Nesterov shows for any initial value $ L_0 > 0 $, 
\begin{equation}  \label{eqn.gg}
  N_k \leq 4(k+1) + 2 \log_2 (L_{k+1}/L_0) \; . 
 \end{equation}

Assume $ L_0 > L_f $.  Since $ k' $ is the first index $ k $ for which $ L_k \leq L_f $ (Lemma~\ref{lem.gd}), we know for $ k < k' $ that $ L_{k+1} = \smfrac{1}{2} L_k $, and hence,
\begin{equation}  \label{eqn.gh}
  \smfrac{1}{2} L_f < L_{k'} = (\smfrac{1}{2} )^{k'} L_0 \leq L_f \; . 
  \end{equation}
Since  $ L_k/L_{k'} $ is a power of 2 for all $ k $, and since $ L_k \leq L_f $ for all $ k > k' $ (Lemma~\ref{lem.gd}), it follows from (\ref{eqn.gh})  that 
\begin{equation} \label{eqn.gi}
 k \geq k' \quad \Rightarrow \quad  L_k \leq L_{k'} \; . 
 \end{equation}

From (\ref{eqn.gg}), (\ref{eqn.gh})  and (\ref{eqn.gi})  we find
\begin{align*}   
     k \geq k'  \quad \Rightarrow \quad      N_k & \leq 4( k + 1 ) + 2 \log_2 ( L_{k+1}/ L_{k'} ) + 2 \log_2 ( L_{k'}/L_0)  \\  & \leq 4(k+1) + 0 - 2k'  \\
                                          & = 4(k - k' + 1) + 2 k'  \\
                                          & = 4(k - k' + 1) +  2 \lceil \log_2(L_0/L_f) \rceil \; ,
\end{align*}
completing the proof.  \hfill $ \Box $

\bibliographystyle{plain}
\bibliography{accelerated_methods}

\end{document}